\documentclass[12pt,a4paper]{article}
\usepackage[dvipdfmx]{graphicx}
\usepackage{latexsym}
\usepackage{amssymb}
\usepackage{amsmath}
\usepackage{amscd}
\usepackage{amsthm}
\usepackage{amsfonts}
\usepackage{mathrsfs}
\usepackage{enumerate}
\usepackage{bm}
\usepackage[usenames]{color}
 \makeatletter
    
    \@addtoreset{equation}{section}
  \makeatother
\newtheorem{Th}{Theorem}[section]

\newtheorem{Pro}[Th]{Proposition}
\theoremstyle{definition}
\newtheorem{Exa}[Th]{Example}
\newtheorem{Rem}[Th]{Remark}
\newcommand{\demo}{\par\noindent{\it Proof. \/}\ }
\newcommand{\enD}{\hfill $\Box$\vspace{3truemm} \par}

\newcommand{\R}{{\mathbb R}}
\newcommand{\lon}{\longrightarrow}

\begin{document}
\title{Pedals and inversions of quadratic curves
 }
\author{Shyuichi IZUMIYA and Nobuko TAKEUCHI }

\date{\today}
\maketitle
\begin{abstract}
The pedal of a curve in the Euclidean plane is a classical subject
which has a singular point at the inflection point of the original curve
or the pedal point.
The primitive of a curve is a curve given by the inverse construction for
making the pedal. In this paper we consider the pedal of a quadratic curve.
On of the main results gives a characterization of such curves, which is 
one of the generalizations of
lima\c{c}ons of Pascal.
\end{abstract}
\renewcommand{\thefootnote}{\fnsymbol{footnote}}
\footnote[0]{2010 Mathematics Subject classification. Primary 53A04;
Secondary 14H50 } \footnote[0]{Key Words and Phrases. plane curves, quadratic curves, quartic curves,
primitive, pedal, inversion} 

\section{Introduction} 
For a regular curve in the Euclidean plane, the {\it pedal} of the curve with respect to the origin (simply, we call the pedal) is the locus of the points on the tangent line of the curve, which are given by the projection image of the origin (the pedal point) along the normal direction.
It might be said the point on the pedal determines the tangent line of the original curve
at the corresponding point.
Therefore, the passage to the pedal is sometimes called the formation of a \lq\lq derivative\rq\rq.
The inverse operation is called the formation of a\lq\lq primitive\rq\rq.
The {\it primitive} of a curve in the plane is the envelope of the normal lines to its 
position vectors at their ends (cf.\cite[pp. 91]{Arnold}).
It is known that  the singularities of the pedal correspond to 
the inflection points of the original curve or the origin.
In \cite{I-T} we introduced the notion of the anti-pedal of a curve
which is the inversion image of the pedal of the curve, so that
singularities of the anti-pedal also correspond to the inflection points of the original curve.
Moreover, we have shown that the primitive is equal to the anti-pedal of the inversion image
of the original curve. 
\par
On the other hand, a {\it lima\c{c}on} (\'Etienne Pascal, 1630) is known to be the pedal of 
a circle.
In this paper we consider a slightly generalized notion of lima\c{c}ons which are the pedals of quadratic curves.
We can show that the pedal and  the inversion image of an irreducible quadratic curve are
cubic or quartic curves.
Since an irreducible quadratic curve is non-singular, the pedal is always well-defined.
However, cubic or quartic curves usually have singularities, so that we need to extend the
definition of pedals and primitives to a class of curves with singularities.
The notion of frontals is one of the candidates of a class of singular curves in the Euclidean plane such that we can develop the differential geometry \cite{F-T, F-T2}.
In \cite{Izu-Take} we adopted the notion of frontals for defining the pedals, the primitives and
the anti-pedals of singular curves.
It is known that plane algebraic curves are frontal (cf. \S 3), so that  the pedals and the inversion
images of irreducible quadratic curves are frontals.
Moreover, we show that the anti-pedal of an irreducible quadratic curve is an irreducible quadratic curve
(cf. Proposition 3.4).
Main results in this paper are Theorems 3.6 and 3.7, which give a characterization of pedals of irreducible quadratic curves. As a consequence, the pedals of irreducible quadratic curves are characterized by
the classification of the anti-pedal as irreducible quadratic curves.
If we start a circle, then we have the lima\c{c}on. Therefore, we call the pedal of an irreducible
quadratic curve a {\it lima\c{c}onoid}.
It has taken around 400 years since lima\c{c}ons were systematically investigated
by \'Etienne Pascal. We can not find such a kind of generalizations so far as the authors know.
\par
We assume that all parametrized curves are class $C^\infty$ unless otherwise stated.

\section{Pedals, primitives and inversions}
In this section we briefly review the properties of pedals and primitives.
Let $\bm{\gamma}:I\lon \R^2$ be a unit speed plane curve, where $\R^2$ is the Euclidean plane with the
canonical scaler product $\langle \bm{a},\bm{b}\rangle=a_1b_1+a_2b_2$ for $\bm{a}=(a_1,a_2),\bm{b}=(b_1,b_2)\in \R^2.$
Then we have the Frenet formulae:
\[
\left\{
\begin{array}{ll}
\bm{t}'(s)=\kappa(s)\bm{n}(s), \\
\bm{n}'(s)=-\kappa (s)\bm{t}(s),
\end{array}
\right.
\]
where $\bm{t}(s)=\bm{\gamma}'(s)$ is the {\it unit tangent vector}, $\bm{n}(s)=J\bm{t}(s)$ is the {\it unit normal vector} and
$\kappa (s)=x_1'(s)x_2''(s)-x_1''(s)x_2'(s)$ is the {\it curvature} of $\bm{\gamma}(s)$, where  $J= \begin{pmatrix} 0 & -1 \\ 1 & 0 \end{pmatrix}$ and $\bm{\gamma}(s)=(x_1(s),x_2(s)).$
The pedal curve of $\bm{\gamma}$ is defined to be
${\rm Pe}_{\bm{\gamma}}(s)=\langle \bm{\gamma}(s),\bm{n}(s)\rangle \bm{n}(s)$ (cf. \cite[Page 36]{Bru-Gib}).
Since
\[
{\rm Pe}_{\bm{\gamma}}'(s)=-\kappa (s)(\langle \bm{\gamma} (s),\bm{t}(s)\rangle\bm{n}(s)+\langle \bm{\gamma}(s),\bm{n}(s)\rangle\bm{t}(s)),
\]
the singular points of the pedal of $\bm{\gamma}$ is the point $s_0$ where  $\bm{\gamma}(s_0)=\bm{0}$ or
$\kappa (s_0)=0$ (i.e. the inflection point).
Therefore, we have the following proposition.
\begin{Pro}
If $\bm{\gamma}$ does not have inflection point, the singular points of the pedal 
${\rm Pe}_{\bm{\gamma}}$ are the points of $\bm{\gamma}$ passing through the origin.
\end{Pro}
By definition, ${\rm Pe}_{\bm{\gamma}}(s)$ is the point on the tangent line through $\bm{\gamma}(s)$, which is given by the projection image of $\bm{\gamma}(s)$ of the normal direction. 
Thus, ${\rm Pe}_{\bm{\gamma}}(s)-\bm{\gamma}(s)$ generates the tangent line at $\bm{\gamma}(s).$
The primitive of a curve $\bm{\gamma}$ in the plane is the envelope of the normal lines to its 
position vectors $\bm{\gamma}(s)$ at their ends (cf.\cite[p. 91]{Arnold}).
It is known that the primitive ${\rm Pr}_{\bm{\gamma}}:I\lon \R^2\setminus \{\bm{0}\}$ of 
$\bm{\gamma}$ with $\langle \bm{\gamma}(s),\bm{n}(s)\rangle\not= 0$ is given by
\[
{\rm Pr}_{\bm{\gamma}}(s)=2\bm{\gamma}(s)-\frac{\|\bm{\gamma}(s)\|^2}{\langle\bm{n}(s),\bm{\gamma}(s)\rangle}\bm{n}(s).
\]
We define the {\it inversion} $\Psi :\R^2\setminus \{\bm{0}\}\lon \R^2\setminus \{\bm{0}\}$ at the origin with respect to the unit circle by
$\displaystyle{\Psi (\bm{x})=\frac{\bm{x}}{\|\bm{x}\|^2}}.$ 
For a subset $A\subset \R^2,$ we define 
\[
\mathscr{I}(A)=\left\{
\begin{matrix} \Psi(A\setminus \{(0,0)\})\cup \{(0,0)\} & \mbox{if}\ \{(0,0)\}\in A, \\
 \Psi (A) & \mbox{if}\  \{(0,0)\}\notin A. 
 \end{matrix}
 \right.
 \]
 Then we have $\mathscr{I}(\mathscr{I}(A))=A.$
We call $\mathscr{I}(A)$ the {\it inversion} of $A.$
For a curve $\bm{\gamma}$ with $\langle \bm{\gamma}(s),\bm{n}(s)\rangle\not= 0,$
we define a mapping ${\rm APe}_{\bm{\gamma}}:I\lon \R^2\setminus \{\bm{0}\}$ by
\[
{\rm APe}_{\bm{\gamma}}(s)=\frac{1}{\langle \bm{\gamma}(s),\bm{n}(s)\rangle}\bm{n}(s),
\]
which is called an {\it anti-pedal curve} of $\bm{\gamma}.$
By a straightforward calculation, we have 
$\Psi \circ {\rm APe}_{\bm{\gamma}}={\rm Pe}_{\bm{\gamma}}$ and $\Psi \circ {\rm Pe}_{\bm{\gamma}}={\rm APe}_{\bm{\gamma}}.$
In \cite{Izu-Take} we have shown the following proposition.
\begin{Pro} For any unit speed plane curve $\bm{\gamma}:I\lon \R^2$
with $\langle \bm{\gamma}(s),\bm{n}(s)\rangle\not= 0,$
we have
\[
{\rm Pr}_{\bm{\gamma}}(s)={\rm APe}_{\Psi\circ\bm{\gamma}}(s)\ \mbox{and}\ {\rm Pr}_{\Psi\circ\bm{\gamma}}(s)={\rm APe}_{\bm{\gamma}}(s).
\]
\end{Pro}

\par 
Although the pedal of a regular curve generally has singularities, there is a tangent line
at any point, so that the primitive of the pedal is well defined.
We can also define the pedal of the primitive of a regular curve.
By definition, we have ${\rm Pr}_{{\rm Pe}_{\bm{\gamma}}}(s)={\rm Pe}_{{\rm Pr}_{\bm{\gamma}}}(s)=\bm{\gamma}(s).$
These arguments can be naturally interpreted by using the notion of frontals
which are natural singular curves in the Euclidean plane such that we can develop the differential geometry \cite{F-T, F-T2}.
We say that $(\bm{\gamma},\bm{\nu}):I\lon \R^{2}\times S^{1}$ is a {\it Legendrian curve}
if $(\bm{\gamma},\bm{\nu})^{*}\theta =0$, where $\theta$ is the canonical contact
$1$-form on the unit tangent bundle $T_{1}\R^{2}=\R^{2}\times S^{1}$ (cf. \cite{Arnold2}). 
The last condition is equivalent to $\langle \dot{\bm{\gamma}}(t),\bm{\nu}(t)\rangle=0$
for any $t\in I.$
We say that $\bm{\gamma}:I\lon \R^{2}$ is a {\it frontal parametrized curve} (briefly, a {\it frontal}\/) if
there exists $\bm{\nu}:I\lon S^{1}$ such that $(\bm{\gamma},\bm{\nu})$ is a Legendrian curve.
If $(\bm{\gamma},\bm{\nu})$ is an immersion (i.e. a regular Legendrian curve), $\bm{\gamma}$ is said to be a {\it front}.
A differential geometry on frontals was constructed in \cite{F-T2}.
For a Legendrian curve $(\bm{\gamma},\bm{\nu}):I\lon \R^{2}\times S^{1},$
we define a unit vector field $\bm{\mu}(t)=J(\bm{\nu}(t))$ along $\bm{\gamma}.$ 
Then we have the following Frenet type formulae \cite{F-T2}:
\[
\left\{
\begin{array}{ll}
\dot{\bm{\nu}}(t)=\ell(t)\bm{\mu}(t), \\
\dot{\bm{\mu}}(t)=-\ell (t)\bm{\nu}(t),
\end{array}
\right.
\]
where $\ell (t)=\langle \dot{\bm{\nu}}(t),\bm{\mu}(t)\rangle.$
Moreover, there exists $\beta (t)$ such that $\dot{\bm{\gamma}}(t)=\beta(t)\bm{\mu}(t)$
for any $t\in I.$
The pair $(\ell,\beta)$ is called a {\it curvature of the Legendrian curve} $(\bm{\gamma},\bm{\nu}).$
By definition, $t_{0}\in I$ is a singular point of $\bm{\gamma}$ if and only if 
$\beta (t_{0})=0.$ Moreover, for a regular curve $\bm{\gamma}$, $\bm{\mu}(t)=\bm{t}(t)$ and 
$\ell (t)=\|\dot{\bm{\gamma}}(t)\|\kappa (t).$
The Legendrian curve $(\bm{\gamma},\bm{\nu})$ is immersive (i.e. $\bm{\gamma}$ is a front) if and only if 
$(\ell(t),\beta(t))\not= (0,0)$ for any $t\in I .$ 
This also means that $\dot{\bm{\gamma}} (t)\not =0$ or $|\bm{\mu}(t),\dot{\bm{\mu}}(t)|\not=0,$
where $|\bm{x},\bm{y}|=x_1y_2-x_2y_1$ for
$\bm{x}=(x_1.x_2),\bm{y}=(y_1,y_2).$
So the inflection point $t_{0}\in I$ of the frontal $\bm{\gamma}$ is a point $\ell(t_{0})=0$
(i.e. $|\bm{\mu}(t_0),\dot{\bm{\mu}}(t_0)|=0,$).
We say that a singular point $t=t_0$ of $\bm{\gamma}$ is an {\it ordinary cusp} if
the image $\bm{\gamma}(I)$ around $\bm{\gamma}(t_0)$ is locally diffeomorphic to
$\{ (t^2,t^3)\in \R^2\ |\ t\in \R\}$ around the origin.  
It is known that $t=t_0$ is the ordinary cusp if and only if
$\dot{\bm{\gamma}}(t_0)=0,$ $\ddot{\bm{\gamma}}(t_0)\not=0$ and
$|\ddot{\bm{\gamma}}(t_0),\dddot{\bm{\gamma}}(t_0)|\not= 0$ (cf. \cite[p. 101]{Gibson}).
Therefore,  $\bm{\gamma}$ is a front if and only if $\bm{\gamma}$ has only ordinary cusps as singular points.
For more detailed properties of Legendrian curves, see \cite{F-T, F-T2}.
We say that $C\subset \R^2$ is a {\it frontal set} (respectively, a {\it front set}) if
for any $\bm{x}_0\in C,$ 
there exist an open neighborhood $U$ of $\bm{x}_0$, $\varepsilon_i >0$ and a frontals (respectively, a fronts) $\bm{\gamma}_i:(t_i-\varepsilon_i, t_i+\varepsilon_i)\lon \R^2$ $(i=1,\dots ,r)$
such that $\bm{\gamma}_i(t_i)=\bm{x}_0$ and 
\[
\cup_{i=1}^r \bm{\gamma}_i(t_i-\varepsilon_i, t_i+\varepsilon_i)=U\cap C.
\]
\par
In \cite{I-T,LP} the {\it pedal} of a frontal $\bm{\gamma}$ is defined by
\[
\mathcal{P}e_{\bm{\gamma}}(t)=\langle \bm{\gamma}(t), \bm{\nu}(t)\rangle \bm{\nu}(t).
\]
We have shown that if there exist $\delta (t)$ and 
$\bm{\sigma}:I\lon S^{1}$such that $\bm{\gamma} (t)=\delta (t)\bm{\sigma} (t)$
for any $t\in I,$  the pedal $\mathcal{P}e_{\bm{\gamma}}(t)$ of
$\bm{\gamma}$ is a frontal.
Moreover, the {\it anti-pedal} of $\bm{\gamma}$ is defined to be
\[
\mathcal{AP}e_{\bm{\gamma}}(t)=\frac{1}{\langle \bm{\gamma}(t), \bm{\nu}(t)\rangle} \bm{\nu}(t)
\]
under the condition that $\langle \bm{\gamma}(t),\bm{\nu}(t)\rangle\not= 0.$
It is easy to show that $\Psi\circ \mathcal{P}e_{\bm{\gamma}}=\mathcal{AP}e_{\bm{\gamma}}$
and $\Psi\circ \mathcal{AP}e_{\bm{\gamma}}=\mathcal{P}e_{\bm{\gamma}}.$
On the other hand, for a frontal $\bm{\gamma}$
with $\langle
\bm{\gamma}(t),\bm{\nu}(t)\rangle\not=0,$ we define the {\it primitive} of $\bm{\gamma}$
by
\[
\mathcal{P}r_{\bm{\gamma}}(t)=2\bm{\gamma}(t)-\frac{\|\bm{\gamma}(t)\|^2}{\langle
\bm{\gamma}(t),\bm{\nu}(t)\rangle}\bm{\nu}(t).
\]
We say that the frontal {\it $\bm{\gamma}$ satisfies the condition $(^*)$}
if $\langle
\bm{\gamma}(t),\bm{\nu}(t)\rangle\not=0.$
We remark that the condition $(^*)$ is called a {\it no-silhouette condition} in
\cite{KN}.
Then we also have shown the following proposition in \cite[Lemmas 5.1 and 5.2]{Izu-Take}.
\begin{Pro} Let $(\bm{\gamma},\bm{\nu}):I\lon \R^{2}\times S^{1}$ be a Legendrian curve
with the condition $(^*).$ Then we have
\[
\mathcal{P}r_{\bm{\gamma}}(t)=\mathcal{AP}e_{\Psi\circ\bm{\gamma}}(t)=\Psi\circ\mathcal{P}e_{\Psi\circ\bm{\gamma}}(t)\ ({\it i.e.}\ \mathcal{P}r_{\Psi\circ \bm{\gamma}}(t)=\mathcal{AP}e_{\bm{\gamma}}(t)=\Psi\circ\mathcal{P}e_{\bm{\gamma}}(t) ).
\]

Moreover, the primitive $\mathcal{P}r_{\bm{\gamma}}$ of 
$\bm{\gamma}$ is a frontal.
\end{Pro}

\section{Quadratic curves}
In this section we consider a general quadratic curve defined by
\[
C:\ g(x,y)=a_{11}x^2+a_{22}y^2+2a_{12}xy+2a_1x+2a_2y+c=0.
\]
Associated with the above quadratic curve, we have two basic invariants:
\[
\Delta_0=\begin{vmatrix} a_{11} & a_{12} \\ a_{12} & a_{22}\end{vmatrix},\ 
\Delta =\begin{vmatrix}
a_{11} & a_{12} & a_1 \\
a_{12} & a_{22} & a_2 \\
a_1 & a_2 & c
\end{vmatrix}.
\]
We say that $C$ is {\it reducible} if there exist lines $\ell_1,\ell _2$ (possibly, $\ell_1=\ell_2$) such that $C=\ell_1\cup\ell_2.$ Then $C$ is said to be {\it irreducible} if it is not reducible.
It is known that $C$ is irreducible if and only if $\Delta \not= 0.$
We suppose that $C$ is irreducible. By the classification theorem on quadratic curves, we have
the following well-known theorem (cf. \cite[p.316]{P-M}).
\begin{Th}
Let $C$ be an irreducible quadratic curve defined by $g(x,y)=0.$ Then
\par\noindent
{\rm (1)} $C$ is an ellipse if $\Delta_0>0$ and $a_{11}\Delta <0,$
\par\noindent
{\rm (2)} $C$ is a hyperbola if $\Delta _0<0$ and $\Delta\not= 0,$
\par\noindent
{\rm (3)} $C$ is a parabola if $\Delta _0=0$ and $\Delta \not=0.$
\end{Th}
For any $(x_0,y_0)\in C,$ the tangent line of $C$ at $(x_0,y_0)$ is
\[
\ell :\  (a_1+a_{11}x_0+a_{12}y_0)x+(a_2+a_{12}x_0+a_{22}y_0)y+c+a_1x_0+a_2y_0=0.
\]
The line through the origin which is orthogonal to $\ell$ is
parametrized by
\[
L:\ \begin{pmatrix} x \\ y \end{pmatrix}
=
\begin{pmatrix}
t(a_1+a_{11}x_0+a_{12}y_0) \\
t(a_2+a_{12}x_0+a_{22}y_0)
\end{pmatrix}.
\]
Then the intersection $\ell\cap L$ is 
\[
\begin{pmatrix}
x \\
y
\end{pmatrix}
=
\begin{pmatrix}
\displaystyle{(a_1+a_{11}x_0+a_{12}y_0)\frac{-c-a_1x_0-a_2y_0}{(a_1+a_{11}x_0+a_{12}y_0)^2+(a_2+a_{12}x_0+a_{22}y_0)^2}}\\
\displaystyle{(a_2+a_{12}x_0+a_{22}y_0)\frac{-c-a_1x_0-a_2y_0}{(a_1+a_{11}x_0+a_{12}y_0)^2+(a_2+a_{12}x_0+a_{22}y_0)^2}}
\end{pmatrix}
\]
and
$
a_{11}x_0^2+a_{22}y_0^2+2a_{12}x_0y_0+2a_1x_0+2a_2y_0+c=0.
$
By definition, the locus of the above points $(x,y)$ is the pedal of $C.$
We now define a polynomial $G_1[C](x,y)$ of degree $4$
by
\begin{eqnarray*}
&{}&G_1[C](x,y)=(a_{12}^2-a_{11}a_{22})(x^2+y^2)^2+2(a_{12}a_2-a_1a_{22})(x^2+y^2)x \\
&{}&\quad +2(a_1a_{12}-a_2a_{11})(x^2+y^2)y 
+(a_2^2-a_{22}c)x^2+2(a_{12}c-a_1a_2)xy+(a_1^2-a_{11}c)y^2.
\end{eqnarray*}
\begin{Pro}
With the same notations as above, 
\[
{\rm Pe}_C=\{(x,y)\in \R^2\ |\ G_1[C](x,y)=0\}.
\]
Therefore ${\rm Pe}_C$ is a cubic or quartic curve through the origin
including the case when the origin is the isolated points on ${\rm Pe}_C.$
\end{Pro}
\demo
We substitute $x_0=X,y_0=Y$ in the above relations.
Then, by a straight forward calculation, we have
\begin{eqnarray*}
x&=&\frac{-a_2^2X+a_{22}X(c+a_1X)+a_1a_2Y-a_{12}cY+a_1a_{22}Y^2-a_{11}a_2(X^2+Y^2)}{-a_1a_{22}X+
a_{12}(a_2X+a_1Y)+a_{12}^2(X^2+Y^2)-a_{11}(a_2Y+a_{22}(X^2+Y^2))}, \\
y&=& \frac{-a_1a_2X+a_{12}cX+a_1a_{12}X^2-a_{11}a_2X^2+a_1^2Y-a_{11}cY+a_1a_{12}Y^2-a_{11}a_2Y^2}{
-a_{12}a_2X+a_1a_{22}X-a_{12}^2X^2+a_{11}a_{22}X^2-a_1a_{12}Y+a_{11}a_2Y-a_{12}^2Y^2+a_{11}a_{22}Y^2}.
\end{eqnarray*} 
We also substitute the above equalities to $g(x,y)=0.$ Then we obtain 
$G_1[C](X,Y)=0.$ 
By the definition of $G_1[C](x,y),$ we have $G_1[C](0,0)=0.$
If $\Delta_0=-(a_{12}^2-a_{11}a_{22})=0,$ then the degree of $G_1[C]$ is three.
 Here, the $(i,j)$ cofactor of $\Delta$ is denoted by $D_{ij}.$ Then we have
\[
\Delta =a_1D_{31}+a_2D_{32}+cD_{33}\ \mbox{and}\ \Delta_0=D_{33}.
\]
Moreover, we have
\[
G_1[C](x,y)=-D_{33}(x^2+y^2)^2+2D_{31}(x^2+y^2)x +2D_{32}(x^2+y^2)y 
-D_{11}x^2-2D_{12}xy-D_{22}y^2.
\]
If the degree of $G_1[C]$ is two, then we have $D_{31}=D_{32}=D_{33}=0,$ so that
$\Delta=0.$ This contradicts to the assumption that $C$ is irreducible.
\enD

We now consider the anti-pedal ${\rm APe}_C$ of C (i.e, the inversion image of ${\rm Pe}_C$).
By definition $\Psi (x,y)=(x/(x^2+y^2),y/(x^2+y^2))=(X,Y),$
so that we have $(x,y)=(X/(X^2+Y^2),X/(X^2+Y^2)$ for $(x,y)\not =(0,0),(X,Y)\not= 0.$
Then we have
\begin{eqnarray*}
&{}& G_1[C]\circ \Psi ^{-1}(X,Y)=(a_2^2-a_{22}c)X^2+(a_1^2-a_{11}c)Y^2+2(a_{12}c-a_1a_2)XY \\
&{}& \quad +2(a_{12}a_2-a_1a_{22})X+2(a_1a_{12}-a_{11}a_2)Y+a_{12}^2-a_{11}a_{22}.
\end{eqnarray*}
Therefore, we have the following proposition.
\begin{Pro}
With the same notations as above, we have
\[
\Psi({\rm Pe}_{C}\setminus \{(0,0)\})=\{(x,y)\in \R^2\setminus \{(0,0)\}\ |\ G_2[C](x,y)=0\},
\] where
\begin{eqnarray*}
&{}& G_2[C](x,y)=(a_2^2-a_{22}c)x^2+(a_1^2-a_{11}c)y^2+2(a_{12}c-a_1a_2)xy \\
&{}& \quad +2(a_{12}a_2-a_1a_{22})x+2(a_1a_{12}-a_{11}a_2)y+a_{12}^2-a_{11}a_{22}.
\end{eqnarray*}
 \end{Pro}
 Since the inversion $\Psi$ is defined on $\R^2\setminus\{(0,0)\}$, the anti-pedal of $C$ is
 defined by $\Psi({\rm Pe}_C\setminus\{(0,0)\}).$
 However, we can generalize the definition of the anti-pedal of the quadratic curve $C$ on the hole plane $\R^2.$
We define ${\rm APe}_{C}=\mathscr{I}({\rm Pe}_C)=\{(x,y)\in \R^2\ |\ G_2[C](x,y)=0\},$ which is called
the {\it anti-pedal} of the quadratic curve $C.$
The basic invariants of the quadratic curve $G_2[C](x,y)=0$ is denoted by $\widehat{\Delta}_0,\widehat{\Delta}.$
By straight forward calculations, we have
\[
\widehat{\Delta}_0=c\Delta,\ \widehat{\Delta}=-\Delta^2.
\]
Thus we have the following proposition.
\begin{Pro} If $C$ is an irreducible quadratic curve, then the anti-pedal ${\rm Ape}_C$ of $C$ is
an irreducible quadratic curve. Moreover, we have the following assertions\/{\rm :}
\par\noindent
{\rm (1)} If the origin is an isolated point of ${\rm Pe}_C,$ then ${\rm Ape}_C$ is an ellipse.
\par\noindent
{\rm (2)}  The anti-pedal ${\rm Ape}_C$ is a parabola if and only if $\{(0,0)\}\in C.$
\end{Pro}
\demo
By the assumption that $C$ is irreducible, $\Delta\not=0,$
so that $\widehat{\Delta}=-\Delta^2\not= 0.$ This means that
${\rm Ape}_C$ is irreducible.
If the origin is an isolated point of the pedal ${\rm Pe}_C$,
the other irreducible components do not pass through the origin.
Therefore, the inversion image ${\rm Ape}_C=\Psi({\rm Pe}_C)$ is a compact irreducible quadratic curve.
This means that it is an ellipse.
Moreover,  $\{(0,0)\}\in C$ if and only if $c=0,$
which means that $\widehat{\Delta}_0=0.$
The last condition is equivalent to that ${\rm Ape}_C$ is a parabola.
This completes the proof. \enD
\par
On the other hand, we consider the singularities of the pedal of a quadratic curve.
Since an irreducible quadratic curve has no inflection points, a singularity of the pedal ${\rm Pe}_C$  exists if and only if the curve $C$ passes through the origin (cf. Proposition 2.1).
By Proposition 3.4, $\{(0,0)\}\in C$ if and only if the anti-pedal ${\rm APe}_C$ is a parabola.
Thus we have the following theorem.
\begin{Th} Let $C$ be an irreducible quadratic curve. Then the following conditions are equivalent\/{\rm :}
\par
{\rm (1)} the quadratic curve $C$ passes through the origin,
\par
{\rm (2)} the pedal ${\rm Pe}_C$ of $C$ has a singular point,
\par
{\rm (3)} the anti-pedal ${\rm APe}_C$ is a parabola.
\end{Th}
\par
Moreover, if we substitute the relation  $(x,y)=(X/(X^2+Y^2),Y/(X^2+Y^2)$ for
$(X,Y)\not= 0, (x,y)\not= 0$ to $g(x,y)=0,$
then we have
\[
c(X^2+Y^2)^2+2a_1X(X^2+Y^2)+2a_2Y(X^2+Y^2)+a_{11}X^2+2a_{12}XY+a_{22}Y^2=0.
\]
We consider a polynomial of degree four defined by
\[
G_3[C](x,y)=c(x^2+y^2)^2+2a_1x(x^2+y^2)+2a_2y(x^2+y^2)+a_{11}x^2+2a_{12}xy+a_{22}y^2.
\]
Then we have a curve defined by $V(G_3[C])=\{(x,y)\ |\ G_3[C](x,y)=0\}.$
If $c=0,$ then $V(G_3[C])$ is a cubic curve. The condition $c=0$ means that $C$ passes through the origin.
If $c=a_1=a_2=0,$ then $V(G_3[C])$ is a quadratic curve. However, $c=a_1=a_2=0$ implies
$\Delta =0,$ this means that $C$ is not irreducible.
Therefore, $V(G_3[C])$ is a cubic or quartic curve.
By the above arguments, $\Psi(C\setminus \{(0,0)\})\subset V(G_3[C]).$ Since $V(G_3[C])$ is a
cubic or  quartic curve, it has singularities generally.
Since $C$ is an irreducible quadratic curve, there are no singular points on $C.$
Moreover, $\Psi$ is a diffeomorphism, so that $\Psi (C\setminus \{(0,0)\})$ is  non-singular
in $\R^2\setminus \{(0,0)\}.$
Therefore, the singular point of $V(G_3[C])$ is the origin.
Since ${\rm Pr}_{\Psi(C\setminus \{(0,0)\})}={\rm APe}_{C\setminus \{(0,0)\}}$, ${\rm Pr}_{\Psi(C\setminus \{(0,0)\})}\subset V(G_2[C])=\{(x,y)\in \R^2\ |\ G_2[C](x,y)=0\}.$
\par
On the other hand, we need some general properties of algebraic curves.
We consider algebraic curves whose primitives are irreducible quadratic curves.
Let
\[
F(x,y)=\sum_{0\leq i+j\leq n} a_{ij}x^iy^j
\]
 be a real polynomial of degree $n$.
We substitute $x=X/(X^2+Y^2),y=Y/(X^2+Y^2),$ so that
\[
F\left(\frac{X}{X^2+Y^2},\frac{Y}{X^2+Y^2}\right)
=\sum_{0\leq i+j\leq n} a_{ij} \left(\frac{X}{X^2+Y^2}\right)^i\left(\frac{Y}{X^2+Y^2}\right)^j.
\]
Therefore $F\left(\frac{X}{X^2+Y^2},\frac{Y}{X^2+Y^2}\right)=0$ if and only if
\[
\sum_{0\leq i+j\leq n} a_{ij}X^iY^j(X^2+Y^2)^{n-(i+j)}=0.
\]
We set
\[
\widetilde{F}(x,y)=\sum_{0\leq i,j\leq n} a_{ij}x^iy^j(x^2+y^2)^{n-(i+j)}.
\]
For an algebraic curve $\mathscr{C}=\{(x,y)\ |\ F(x,y)=0\},$ the inversion of $\mathscr{C}$
is 
$\mathscr{I}(\mathscr{C})=\{(x,y)\ |\ \widetilde{F}(x,y)=0\}.$
By the converse construction, we have
$\widetilde{\widetilde{F}}(x,y)=F(x,y),$
so that $\mathscr{I}(\mathscr{I}(\mathscr{C}))=\mathscr{C}.$
With this notation, we have $\mathscr{I}(C)=V(G_3[C])$ for an irreducible quadratic curve $C.$
For any non-isolated point $\bm{x}_0$ of the irreducible component of $\mathscr{C},$
there exists $\varepsilon>0$ and an analytic mapping $\bm{\gamma}:(-\varepsilon, \varepsilon)\lon \R^2$
such that $\bm{\gamma}(0)=\bm{x}_0$ (cf. \cite[Lemma 3.3]{Milnor}).
It follows that if an irreducible component of an algebraic curve is non-isolated, then it is a frontal set.
Moreover, an isolated point is a front set by definition.
Therefore, any algebraic curve is a frontal set.
On the other hand, for an algebraic curve $\mathscr{C}$, we define $\mathscr{C}^*$ by
\[
\mathscr{C}^*=\{(x,y)\in \mathscr{C}\ |\ (x,y)\ \mbox{satisfies\ the\ condition}\ (^*)\}.
\]
We remark that $(0,0)\notin \mathscr{C}^*.$ Then $\Psi (\mathscr{C}^*)$ is well-defined.
Since $\Psi$ is a conformal diffeomorphism, $\Psi (\mathscr{C}^*)$ is a frontal and satisfy
the condition $(^*).$
Therefore, $\mathcal{P}r_{\mathscr{C}^*}, \mathcal{AP}e_{\mathscr{C}^*}$ and $\mathcal{AP}e_{\Psi(\mathscr{C}^*)}$ are
well-defined.
However, $\mathscr{I}(C)$ has singularities generally, so that we have to
extend the above notions for singular curves.
In \S 2, we considered the class of frontals in $\R^2$ and generalized all the above notions for
frontals. Here, we show the following theorem.
\begin{Th}
Let $C$ be an irreducible quadratic curve. Then
$\mathscr{I}(C)$ is a cubic or quartic curve and each irreducible component is a front set.
Moreover, we have the following assertions\/{\rm :}
\par\noindent
{\rm (1)} If $C$ is an ellipse, the origin is an isolated point of the inversion $\mathscr{I}(C)$ of $C$ such that $\mathscr{I}(C)\setminus \{(0,0)\}$ is regular as a parametrized curve.
\par\noindent
{\rm (2)}  If $C$ is a hyperbola, the inversion $\mathscr{I}(C)$ of $C$ is regular as a parametrized curve and it has a node at the origin
as a quartic curve.
The tangent directions of the curve at the node are equal to the asymptotic directions of the hyperbola $C$.
\par\noindent
{\rm (3)} If $C$ is a parabola, the inversion $\mathscr{I}(C)$ of $C$ has the ordinary cusp at the origin
such that $\mathscr{I}(C)\setminus \{(0,0)\}$ is regular as a parametrized curve.
The tangent direction of the curve at the cusp is equal to the direction of the axis of the parabola $C$. 
\end{Th}
\demo
By definition $\mathscr{I}(C)$ is a cubic or quartic curve.

Since the notion of frontals is defined for parametrized curve,
we consider a parametrization of $C$.
Suppose that $C$ is a parabola. Then
any parabola can be parametrized by
\[
^t\bm{\gamma}(t)=\begin{pmatrix} \cos\theta & -\sin\theta \\
\sin\theta & \cos\theta 
\end{pmatrix}
\begin{pmatrix} t \\ at^2
\end{pmatrix}
+\begin{pmatrix} b_1 \\ b_2
\end{pmatrix}
\]
for $\theta \in [0,2\pi), a\not=0,b_1,b_2\in \R.$
Thus we have
\[
\bm{\gamma}(t)=(\cos\theta t-a\sin\theta t^2+b_1,\sin\theta t+a\cos\theta t^2+b_2).
\]
Since we need information on points at infinity, we extend the plane $\R^2$ to
the projective plane $P\R^2.$
For the homogeneous coordinates $[X:Y:Z]\in P\R^2,$
we identify the Euclidean plane $\R^2$ to the affine plane
defined by
$Z\not= 0$ (i.e. $\R^2=\{[X:Y:Z]\ |\  Z\not=0\}$)
Then we have
$\displaystyle{(x,y)=\left(\frac{X}{Z},\frac{Y}{Z}\right)}.$
We consider the inversion $\Psi(x,y)=\displaystyle{\left(\frac{x}{x^2+y^2},\frac{y}{x^2+y^2}\right)}.$
Since 
$\displaystyle{
\Psi\left(\frac{X}{Z},\frac{Y}{Z}\right)=\left(\frac{ZX}{X^2+Y^2},\frac{ZY}{X^2+Y^2}\right)
},$
we define
$\widetilde{\Psi}:\R^2\lon P\R^2$
by
\[
\widetilde{\Psi}([X:Y:Z])=[ZX:ZY:X^2+Y^2],
\]
which is a projective version of the inversion.
We substitute $t=\frac{X}{Z}$ in the parabola $\bm{\gamma}$. Then we have
\begin{eqnarray*}
&{}& \left[ \cos\theta \frac{X}{Z}-a\sin\theta \frac{X^2}{Z^2}+b_1:\sin\theta \frac{X}{Z}+a\cos\theta \frac{X^2}{Z^2}+b_2:1\right] \\
&{}& =[\cos\theta ZX-a\sin\theta X^2+b_1Z^2:\sin\theta ZX+a\cos\theta X^2+b_2Z^2:Z^2].
\end{eqnarray*}
We define $\Gamma : P\R^1\lon P\R^2$ by
\[
\Gamma ([X:Z])=[\cos\theta ZX-a\sin\theta X^2+b_1Z^2:\sin\theta ZX+a\cos\theta X^2+b_2Z^2:Z^2],
\]
which gives a parametrization of a parabola in the projective plane $P\R^2.$
Therefore, we have
\[
\widetilde{\Psi}\circ \Gamma ([X:Z])=[\cos\theta Z^3X-a\sin\theta Z^2X^2+b_1Z^4:
\sin\theta Z^3X+a\cos\theta Z^2X^2+b_2Z^4: F(X,Z)],
\]
where
\[
F(X,Z)=(\cos\theta ZX-a\sin\theta X^2+b_1Z^2)^2+(\sin\theta ZX+a\cos\theta X^2+b_2Z^2)^2.
\]
On the other hand, if $X\not=0,$ then for $[X:Z]=[1:Z/X],$ we have
\[
\widetilde{\Psi}\circ \Gamma \left(\left[1:\frac{Z}{X}\right]\right)=
\left[\cos\theta \frac{Z^3}{X^3}-a\sin\theta \frac{Z^2}{X^2}+b_1\frac{Z^4}{X^4}:
\sin\theta \frac{Z^3}{X^3}+a\cos\theta \frac{Z^2}{X^2}+b_2\frac{Z^4}{X^4}:
F\left(1,\frac{Z}{X}\right)\right].
\]
We substitute $\tau =Z/X,$ so that
\[
\widetilde{\Psi}\circ \Gamma \left(\left[1:\tau\right]\right)=
\left[\cos\theta \tau^3-a\sin\theta \tau^2+b_1\tau^4:
\sin\theta \tau^3+a\cos\theta \tau^2+b_2\tau^4:
f\left(\tau\right)\right],
\]
where $f(\tau)=F(1,\tau).$
Since 
\[
f\left(\tau\right)=(\cos\theta \tau-a\sin\theta +b_1\tau^2)^2+(\sin\theta \tau+ a\cos\theta +b_2\tau^2)^2,
\]
 we have $f(0)=a^2\not=0$ and $\widetilde{\Psi}\circ \Gamma \left(\left[1:0\right]\right)=
[0:0:a^2]=[0:0:1].$
For sufficiently small $\varepsilon >0,$ we have $f(\tau)\not= 0$ for $\tau \in (-\varepsilon, \varepsilon).$
Therefore, the inversion of the parabola $\bm{\sigma}=\widetilde{\Psi}\circ \Gamma :(-\varepsilon, \varepsilon)\lon
\R^2\subset P\R^2$ is given by 
\[
\bm{\sigma}(\tau)=\left( \frac{\cos\theta \tau^3-a\sin\theta \tau^2+b_1\tau^4}{f(\tau)},
\frac{\sin\theta \tau^3+a\cos\theta \tau^2+b_2\tau^4}{f(\tau)}\right).
\]
If we define
\begin{eqnarray*}
g_1(\tau)&=&(3\cos\theta \tau -2a\sin\theta+4b_1\tau^2)f(\tau)-(\cos\theta \tau^2-a\sin\theta \tau +b_1\tau^3)
\frac{df}{d\tau}(\tau), \\
g_2(\tau)&=& (3\sin\theta \tau+2a\cos\theta +4b_2\tau^2)f(\tau)-(\sin\theta \tau^2+a\cos\theta \tau +b_2\tau^3)\frac{df}{d\tau}(\tau),
\end{eqnarray*}
then
$\displaystyle{d\frac{\bm{\sigma}}{d\tau} (\tau)=\tau\bm{\mu}(\tau)},$ where
$\displaystyle{
\bm{\mu}(\tau)=\frac{1}{f(\tau)^2}(g_1(\tau),g_2(\tau))}.$
Since
\[
\bm{\mu}(0)=\left(\frac{-2a\sin\theta}{f(0)},\frac{2a\cos\theta}{f(0)}\right)=
\frac{2}{a}(-\sin\theta ,\cos\theta)\not= (0,0),
\]
$\bm{\mu}(\tau)$ gives a tangent direction of $\bm{\sigma}(\tau)$ around $\tau=0.$
This means that $\bm{\sigma}(\tau)$ is a frontal.
Since $\dot{\bm{\sigma}}(\tau)=\tau \bm{\mu}(\tau),$ we have
$\ddot{\bm{\sigma}}(\tau)=\bm{\mu}(\tau)+\tau\dot{\bm{\mu}}(\tau).$
Therefore, $\dot{\bm{\sigma}}(0)=\bm{0}$ and $\ddot{\bm{\sigma}}(0)=\bm{\mu}(0)\not=\bm{0}.$
By straight forward calculations, we have
\[
\dot{g}_1(0)=3\cos\theta f(0)=3a^2\cos\theta\ \dot{g}_2(0)=3\sin\theta f(0)=3a^2\sin\theta
\]
and
\[
\dot{\bm{\mu}}(0)=\frac{1}{f(0)^2}(\dot{g}_1(0),\dot{g}_2(0))=\frac{3}{a^2}(\cos\theta,\sin\theta).
\]
Since $\dddot{\bm{\sigma}}(\tau)=2\dot{\bm{\mu}}(\tau)+\tau \ddot{\bm{\mu}}(\tau),$
we have 
\[
|\ddot{\bm{\sigma}}(0), \dddot{\bm{\sigma}}(0)|=|\bm{\mu}(0),\dot{\bm{\mu}}(0)|
=2\begin{vmatrix} -\frac{2}{a}\sin\theta & \frac{2}{a}\cos\theta \\
\frac{3}{a^2}\cos\theta & \frac{3}{a^2}\sin\theta
\end{vmatrix}=-\frac{12}{a^3}\not= 0.
\]
This means that $\sigma$ is an ordinary cusp at $\tau=0$ such that the tangent direction of the
cusp is the direction of axis of the parabola.
Since $\Psi: \R^2\setminus \{(0,0)\}\lon \R^2\setminus \{(0,0)\}$ is a diffeomorphism
and $C$ is a regular curve, all points of $\Psi (C\setminus \{(0,0)\})\setminus \{(0,0)\}$ are non-singular. Moreover, we cannot define $\Psi$ at the origin. This means that $V(G_3[C])$ has only the ordinary cusp at the origin as the singularities, so that $\Psi(C\setminus \{(0,0)\})$ is a front. 
\par
On the other hand, suppose that $C$ is a hyperbola.
Then we have the normal form of a hyperbola:
\[
\frac{x^2}{a^2}-\frac{y^2}{b^2}=-1, (a,b>0),
\]
then $\displaystyle{y=\pm \frac{b}{a}\sqrt{x^2+a^2}}.$
Thus
any hyperbola can be parametrized by
\[
^t\bm{\gamma}(t)=\begin{pmatrix} \cos\theta & -\sin\theta \\
\sin\theta & \cos\theta 
\end{pmatrix}
\begin{pmatrix} t \\ \displaystyle{\pm \frac{b}{a}\sqrt{t^2+a^2} }
\end{pmatrix}
+\begin{pmatrix} c_1 \\ c_2
\end{pmatrix}
\]
for $\theta \in [0,2\pi), a,b>0,c_1,c_2\in \R.$
Thus we have
\[
\bm{\gamma}(t)=\left(\cos\theta t\mp \sin\theta \frac{b}{a}\sqrt{t^2+a^2}+c_1,\sin\theta t\pm\cos\theta \frac{b}{a}\sqrt{t^2+a^2}+c_2\right).
\]
We also consider in the projective plane $P\R^2.$
Then we substitute $t=\frac{X}{Z}$ to the hyperbola $\bm{\gamma}$. Then we have
\begin{eqnarray*}
&{}& \left[ \cos\theta \frac{X}{Z}\mp \sin\theta \frac{b}{a}\sqrt{\frac{X^2}{Z^2}+a^2}+c_1:\sin\theta \frac{X}{Z}\pm \cos\theta \frac{b}{a}\sqrt{\frac{X^2}{Z^2}+a^2}+c_2:1\right] \\
&{}& =\left[\cos\theta X\mp\sin\theta \frac{b}{a}\sqrt{X^2+a^2Z^2}+c_1Z:\sin\theta X\pm\cos\theta \frac{b}{a}\sqrt{X^2+a^2Z^2}+c_2Z:Z\right].
\end{eqnarray*}
We define $\Gamma : P\R^1\lon P\R^2$ by
\[
\Gamma ([X:Z])=\left[\cos\theta X\mp\sin\theta \frac{b}{a}\sqrt{X^2+a^2Z^2}+c_1Z:\sin\theta X\pm\cos\theta \frac{b}{a}\sqrt{X^2+a^2Z^2}+c_2Z:Z\right],
\]
which gives a parametrization of a hyperbola in the projective plane $P\R^2.$
Thus, $\widetilde{\Psi}\circ \Gamma ([X:Z])$ is equal to 
\[
\left[\cos\theta ZX\mp\sin\theta \frac{b}{a}Z\sqrt{X^2+a^2Z^2}+c_1Z^2:\sin\theta ZX\pm\cos\theta \frac{b}{a}Z\sqrt{X^2+a^2Z^2}+c_2Z^2: F(X,Z)\right],
\]
where
\begin{eqnarray*}
F(X,Z)&=&\left(\cos\theta X\mp\sin\theta \frac{b}{a}\sqrt{X^2+a^2Z^2}
+c_1Z\right)^2 \\ 
&{}& \qquad + 
\left(\sin\theta X\pm\cos\theta \frac{b}{a}\sqrt{X^2+a^2Z^2}+c_2Z\right)^2.
\end{eqnarray*}
On the other hand, if $X\not=0,$ then for $[X:Z]=[1:Z/X],$ then $\widetilde{\Psi}\circ \Gamma \left(\left[1:\frac{Z}{X}\right]\right)$ is
\[
\left[\cos\theta \frac{Z}{X}\mp\sin\theta \frac{b}{a}\frac{Z}{X}\sqrt{1+a^2\frac{Z^2}{X^2}}+c_1\frac{Z^2}{X^2}:
\sin\theta \frac{Z}{X}\pm\frac{b}{a}\frac{Z}{X}\sqrt{1+a^2\frac{Z^2}{X^2}}+c_2\frac{Z^2}{X^2}:
F\left(1,\frac{Z}{X}\right)\right].
\]
We substitute $\tau =Z/X,$ so that
\[
\widetilde{\Psi}\circ \Gamma \left(\left[1:\tau\right]\right)=
\left[\cos\theta \tau\mp \sin\theta\frac{b}{a}\tau\sqrt{1+a^2\tau^2}+c_1\tau^2:
\sin\theta \tau\pm \cos\theta \frac{b}{a}\tau\sqrt{1+a^2\tau^2}+c_2\tau^2:
f\left(\tau\right)\right],
\]
where $f(\tau)=F(1,\tau).$
Since 
\[
f\left(\tau\right)=(\cos\theta \mp \sin\theta\frac{b}{a}\sqrt{1+a^2\tau^2}+c_1\tau)^2+
(\sin\theta \pm \cos\theta\frac{b}{a}\sqrt{1+a^2\tau^2}+c_2\tau)^2,
\]
we can show that  $f(0)\not=0$ and $\widetilde{\Psi}\circ \Gamma \left(\left[1:0\right]\right)=
[0:0:f(0)]=[0:0:1].$
For sufficiently small $\varepsilon >0,$ we have $f(\tau)\not= 0$ for $\tau \in (-\varepsilon, \varepsilon).$
Therefore, the inversion of the hyperbola $\bm{\sigma}=\widetilde{\Psi}\circ \Gamma :(-\varepsilon, \varepsilon)\lon
\R^2\subset P\R^2$ is given by 
\[
\bm{\sigma}(\tau)=\left( \frac{\cos\theta \tau\mp \sin\theta\frac{b}{a}\tau\sqrt{1+a^2\tau^2}+c_1\tau^2}{f(\tau)},
\frac{\sin\theta \tau\pm \cos\theta \frac{b}{a}\tau\sqrt{1+a^2\tau^2}+c_2\tau^2}{f(\tau)}\right).
\]
We can show that
\[
d\frac{\bm{\sigma}}{d\tau} (\tau)=\frac{1}{f(0)^2}\left(\cos\theta \mp\sin\theta \frac{b}{a},\sin\theta\pm \cos\theta \frac{b}{a}\right)\not= (0,0).
\]
Therefore, $\bm{\sigma}(\tau)$ is regular at $\tau=0.$
Since $V(G_3[C])$ is non-singular at $V(G_3[C])\setminus \{(0,0)\},$ $V(G_3[C])$ is a front.
Moreover, $(a\cos\theta\mp\sin\theta b,a\sin\theta \pm \cos\theta b)$ is the asymptotic direction of
the hyperbola $\gamma (t).$
This means that 
the origin is a node of $V(G_3[C])$ such that the tangent directions are
equal to the asymptotic directions of the hyperbola $C.$
\enD
We now apply the above theorem.
\begin{Th}
Let $\mathscr{C}$ be an algebraic curve in $\R^2$ such that $C=\mathcal{P}r_{\mathscr{C}^*}$ is an irreducible quadratic curve.
Then we have the following assertions\/{\rm :}
\par\noindent
{\rm (1)} The algebraic curve $\mathscr{C}$ is a quartic curve or a cubic curve and a front such that 
$\mathscr{I}(\mathscr{C})$ is an irreducible quadratic curve.
Moreover, if the equation of $C$ is 
\[
a_{11}x^2+a_{22}y^2+2a_{12}xy+2a_1x+2a_2y+c=0,
\]
then the equation of $\mathscr{C}$ is 
\begin{eqnarray*}
&{}&(a_{12}^2-a_{11}a_{22})(x^2+y^2)^2+2(a_{12}a_2-a_1a_{22})(x^2+y^2)x \\
&{}&\quad +2(a_1a_{12}-a_2a_{11})(x^2+y^2)y 
+(a_2^2-a_{22}c)x^2+2(a_{12}c-a_1a_2)xy+(a_1^2-a_{11}c)y^2=0.
\end{eqnarray*}
If $C$ is a parabola, then $\mathscr{C}$ is a cubic curve.
\par\noindent
{\rm (2)} If $\mathscr{I}(\mathscr{C})$ is an ellipse, the origin is an isolated point of $\mathscr{C}.$
Moreover, $\mathscr{C}\setminus \{(0,0)\}$ is a regular curve.
\par\noindent
{\rm (3)}  If $\mathscr{I}(\mathscr{C})$  is a hyperbola, $\mathscr{C}$ is regular as a parametrized curve and it has only one node at the origin as a quartic curve.
Moreover, the tangent directions of the curve at the node are equal to the asymptotic directions of the hyperbola $\mathscr{I}(\mathscr{C})$.
\par\noindent
{\rm (4)} If $\mathscr{I}(\mathscr{C})$ is a parabola, $\mathscr{C}$ has the ordinary cusp at the origin
and other part are non-singular.
Moreover, the tangent direction of the curve at the cusp is equal to the direction of the axis of the parabola $\mathscr{I}(\mathscr{C})$. 
\end{Th}
\demo
Since $C=\mathcal{P}r_{\mathscr{C}^*}=\mathcal{AP}e_{\Psi (\mathscr{C}^*)}=\Psi(\mathcal{P}e_{\Psi (\mathscr{C}^*)})$, we have
\[
\mathcal{AP}e_{C}=\mathcal{P}r_{\mathscr{I}(C)}=\mathcal{P}r_{\Psi(\mathcal{P}r_{\mathscr{C}^*})}
 =\mathcal{P}r_{\mathcal{P}e_{\Psi(\mathscr{C}^*)}}=\Psi (\mathscr{C}^*)\subset \mathscr{I}(\mathscr{C}),
 \]
 Since $C$ is an irreducible quadratic curve, Proposition 3.4 asserts that
$\overline{\Psi (\mathscr{C}^*)}=\mathscr{I}(\mathscr{C})$
 is an irreducible quadratic curve.
By Theorem 3.6, $\mathscr{C}=\mathscr{I}(\mathscr{I}(\mathscr{C}))$ is a quartic curve and a front.
 This completes assertion (1).
The other assertions also follow from the fact $\mathscr{C}=\mathscr{I}(\mathscr{I}(\mathscr{C}))$
and Theorem 3.6.
This completes the proof. \enD
\begin{Rem}
In the above theorem, the converse assertion also holds.
For an irreducible quadratic curve $C,$ the pedal $\mathscr{C}=\mathcal{P}e_{C}$ is
the cubic or quartic curve given in the above theorem.
\end{Rem}

\section{Examples}
In this section we give some examples of quartic curves and cubic curves whose primitives are quadratic curves.
\begin{Exa}
We consider an ellipse $C_1$ (Fig. 1) defined by $4x^2+3y^2-2xy-6x-4y-6=0.$
In this case the quartic curve $\mathscr{C}_1$ (Fig. 2) is defined by
\[
(x^2+y^2)^2-2x(x^2+y^2)-2y(x^2+y^2)-2x^2-3y^2=0.
\]
Then the inversion $\mathscr{I}(\mathscr{C}_1)$ (Fig.3) of $\mathscr{C}_1$ is an ellipse
defined by $
2x^2+3y^2+2x+2y-1=0.$
By Theorem 3.7, the quartic curve $\mathscr{C}_1$ contains the origin as an isolated point.
Another component is a regular curve.
The relation of these curves is depicted in Fig. 4.
\begin{center}
\begin{tabular}{cc}   
\includegraphics[width=30mm]{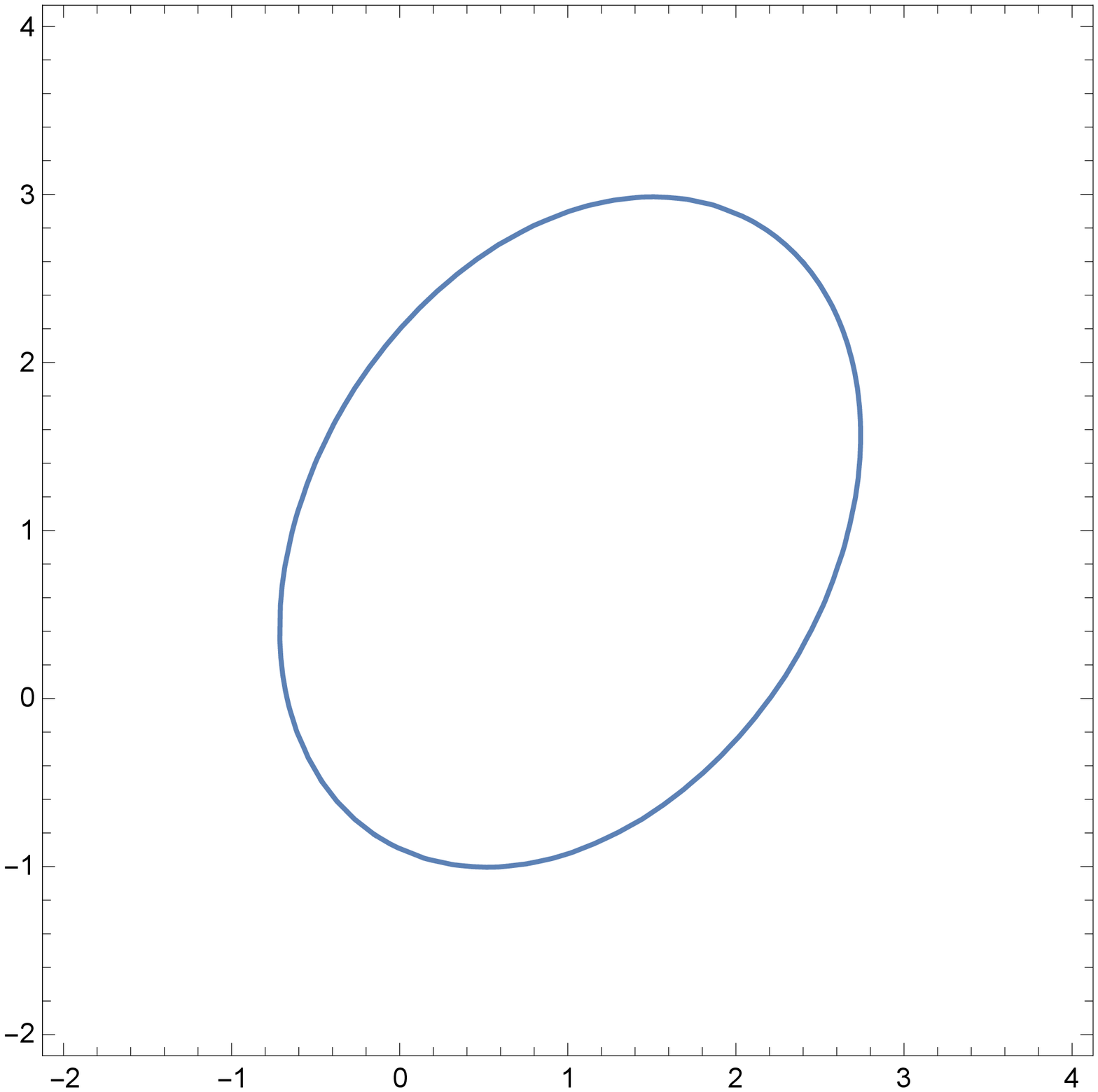}
&
\includegraphics[width=30mm]{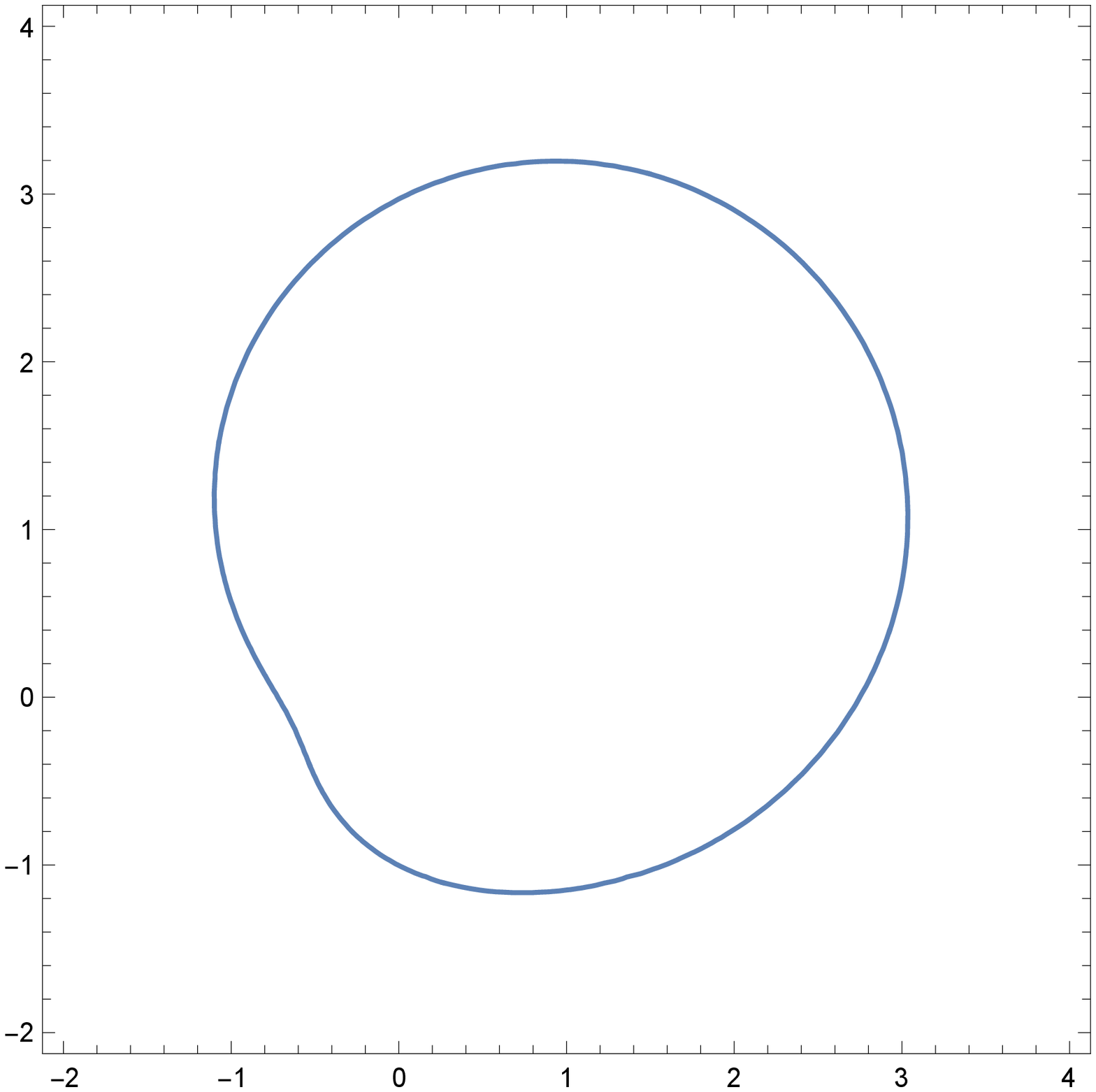}
 \\
\text{Fig. 1 : An ellipse $C_1$}
&
\text{Fig. 2 : The quartic curve $\mathscr{C}_1$}
\end{tabular}
\end{center}

\begin{center}
\begin{tabular}{cc}   
\includegraphics[width=30mm]{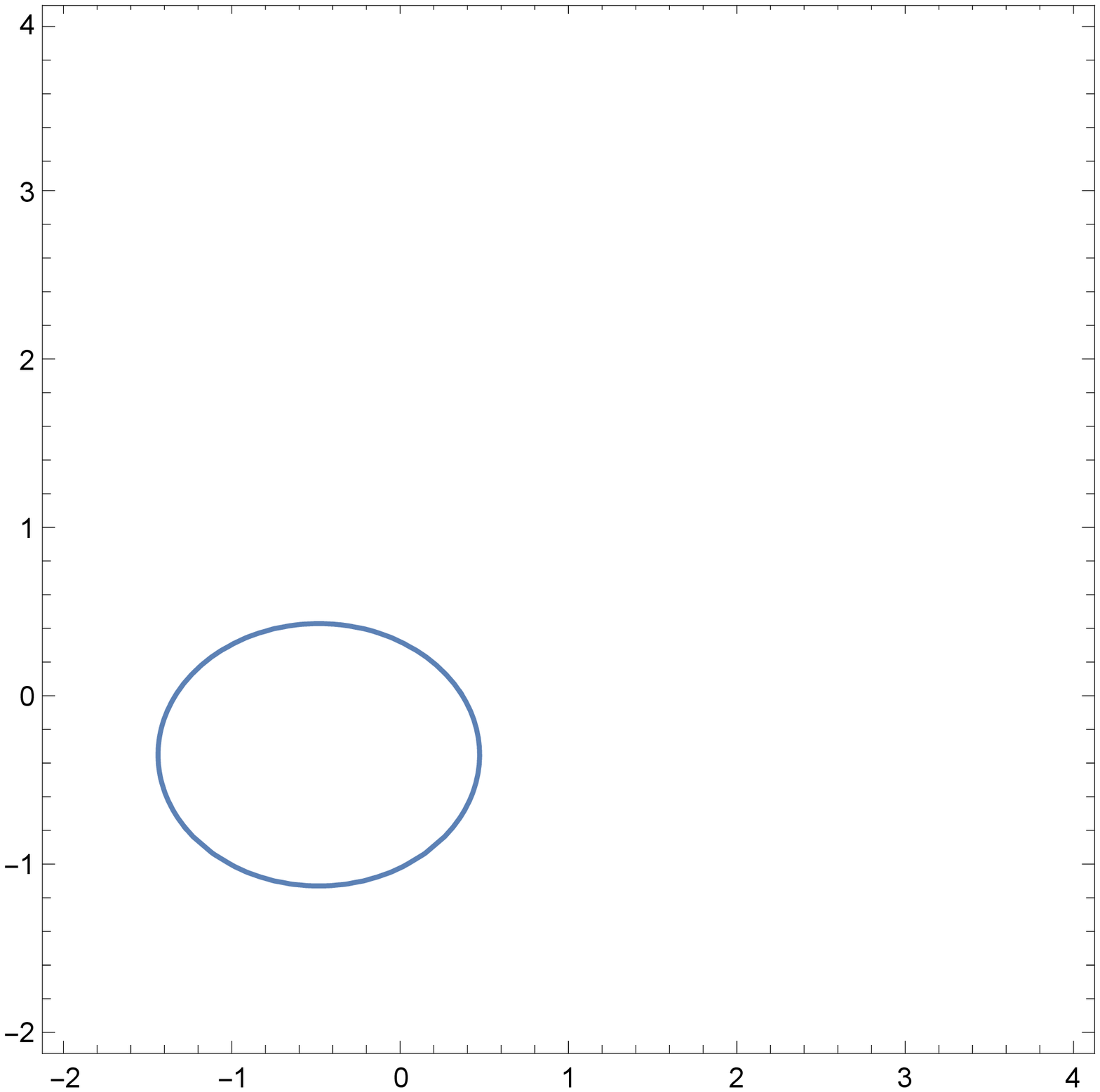}
&
\includegraphics[width=30mm]{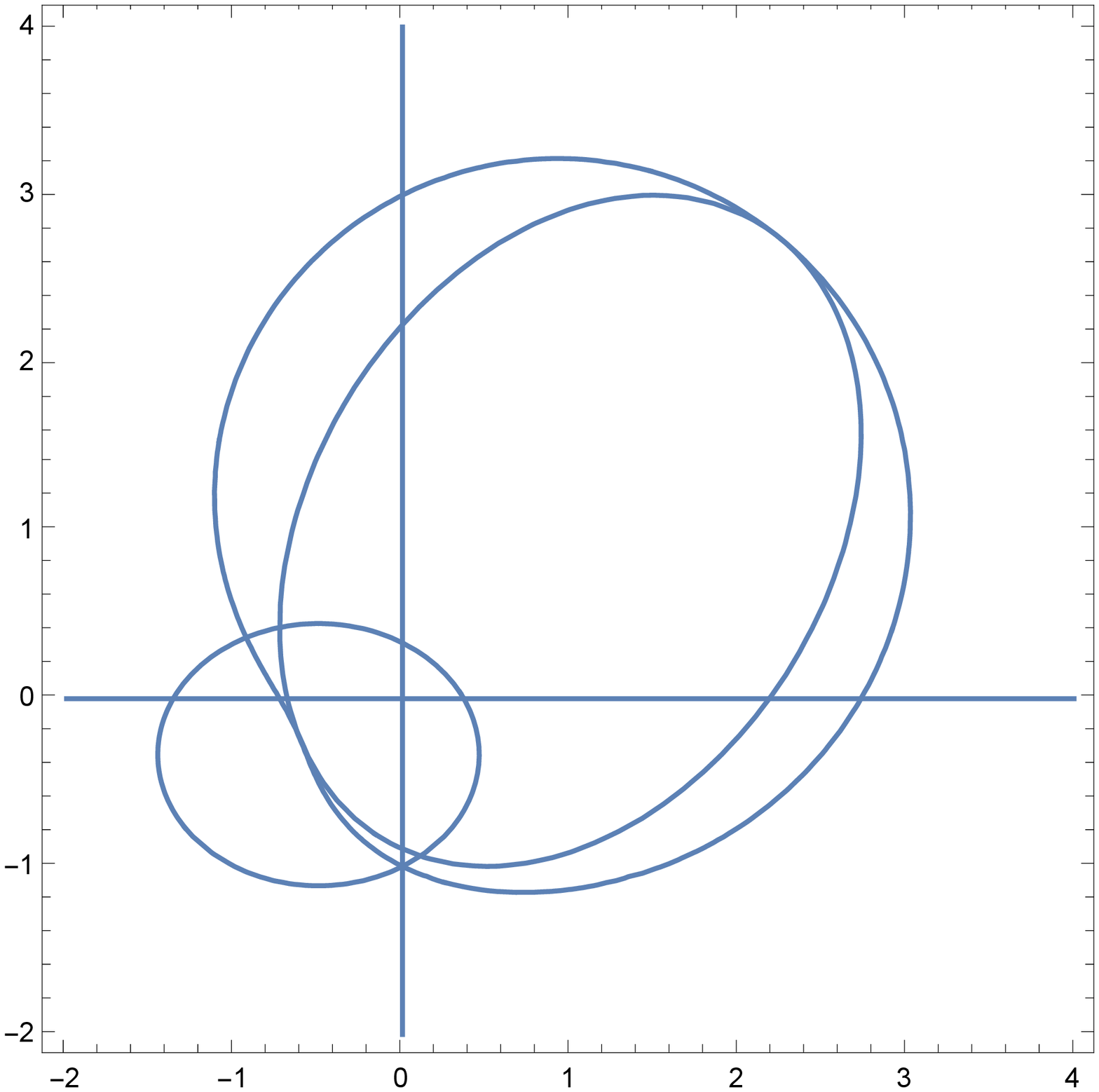}
 \\
\text{Fig. 3 : The ellipse $\mathscr{I}(\mathscr{C}_1)$}
&
\text{Fig. 4 : Altogether}
\end{tabular}
\end{center}

\end{Exa}
\begin{Exa}
We consider a parabola $C_2$ (Fig. 5) defined by $x^2+y^2-2xy+4x+6y+14=0.$
In this case the cubic curve $\mathscr{C}_2$ (Fig. 6) is defined by
\[
2x(x^2+y^2)+2y(x^2+y^2)+x^2+8xy+2y^2=0.
\]
Then the inversion $\mathscr{I}(\mathscr{C}_2)$ (Fig. 7) of $\mathscr{C}_2$ is a hyperbola
defined by $
x^2+2y^2+8xy+2x+2y=0.$
By Theorem 3.7, the cubic curve $\mathscr{C}_2$ has only one node at the origin.
The relation of these curves is depicted in Fig. 8.
\begin{center}
\begin{tabular}{cc}   
\includegraphics[width=30mm]{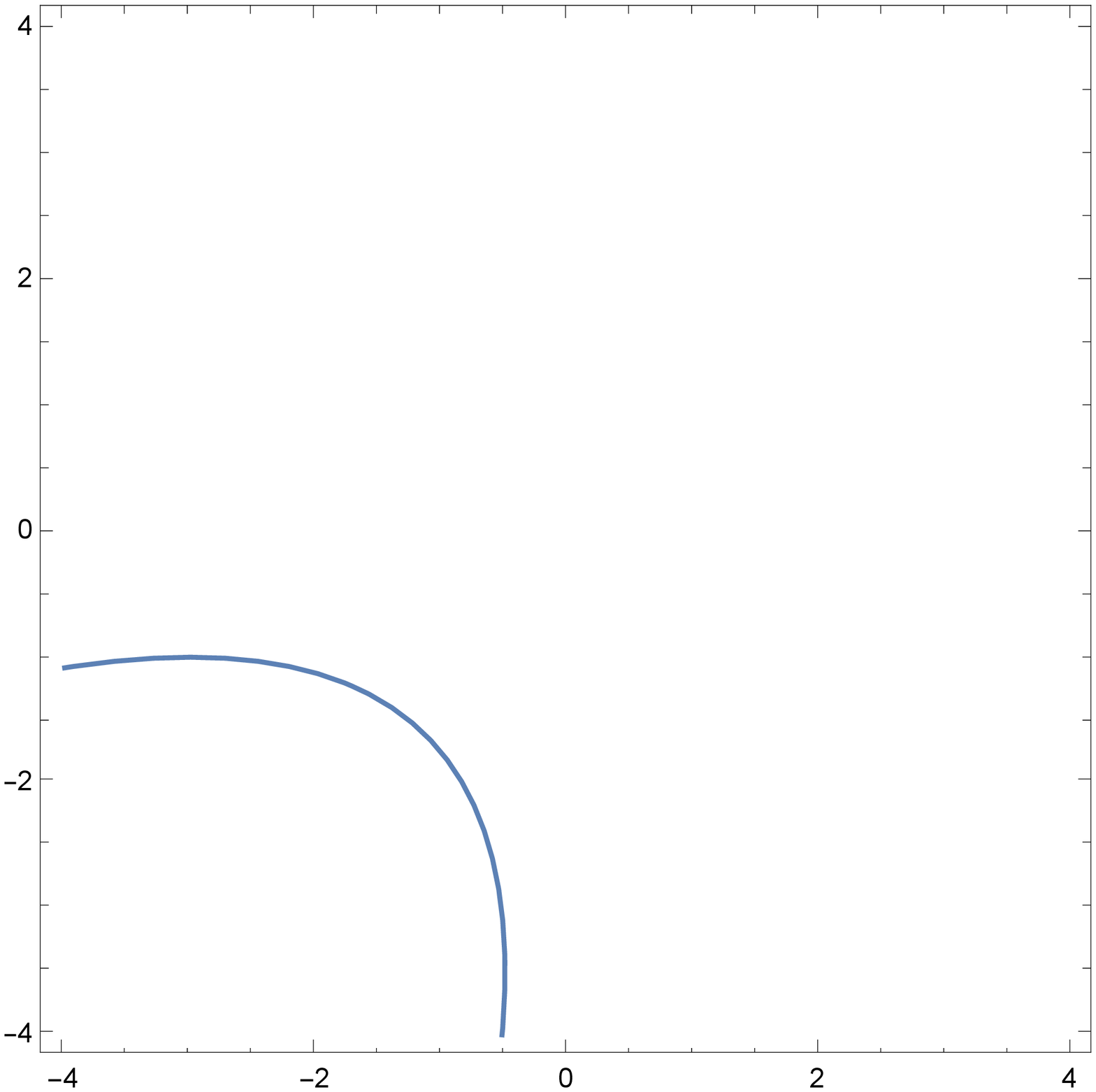}
&
\includegraphics[width=30mm]{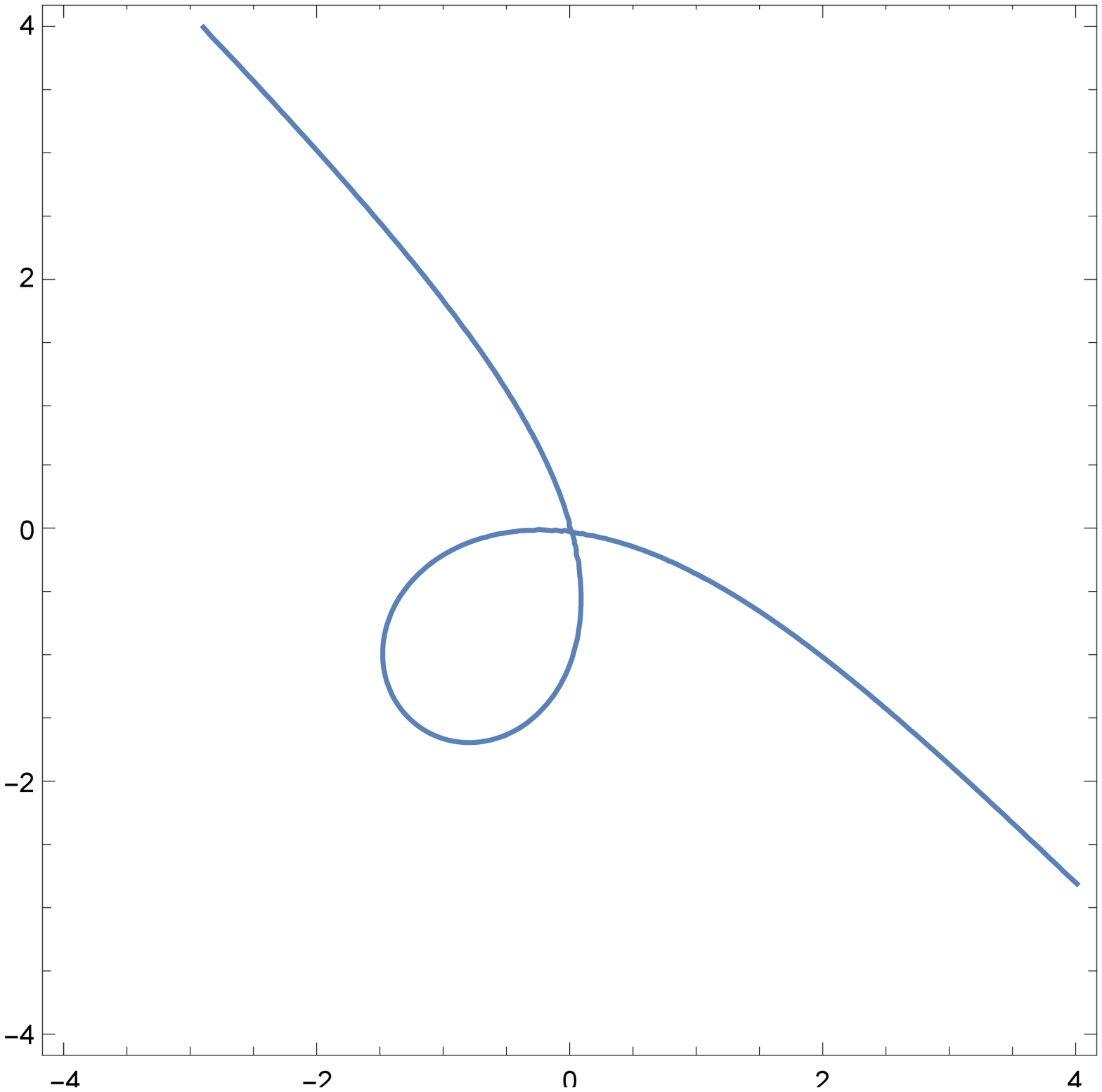}
 \\
\text{Fig. 5 : A parabola $C_2$}
&
\text{Fig. 6 : The cubic curve $\mathscr{C}_2$}
\end{tabular}
\end{center}

\begin{center}
\begin{tabular}{cc}   
\includegraphics[width=30mm]{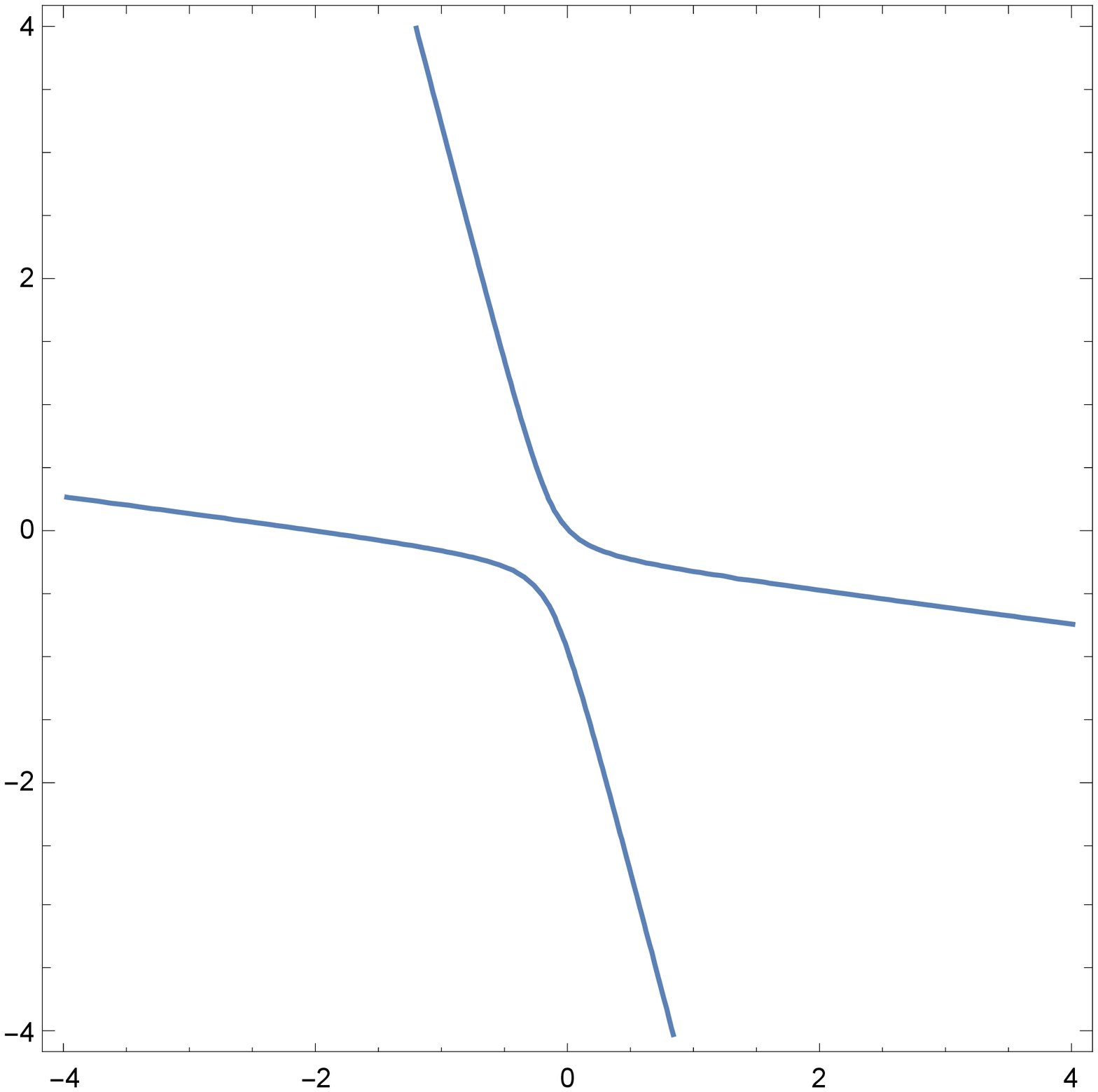}
&
\includegraphics[width=30mm]{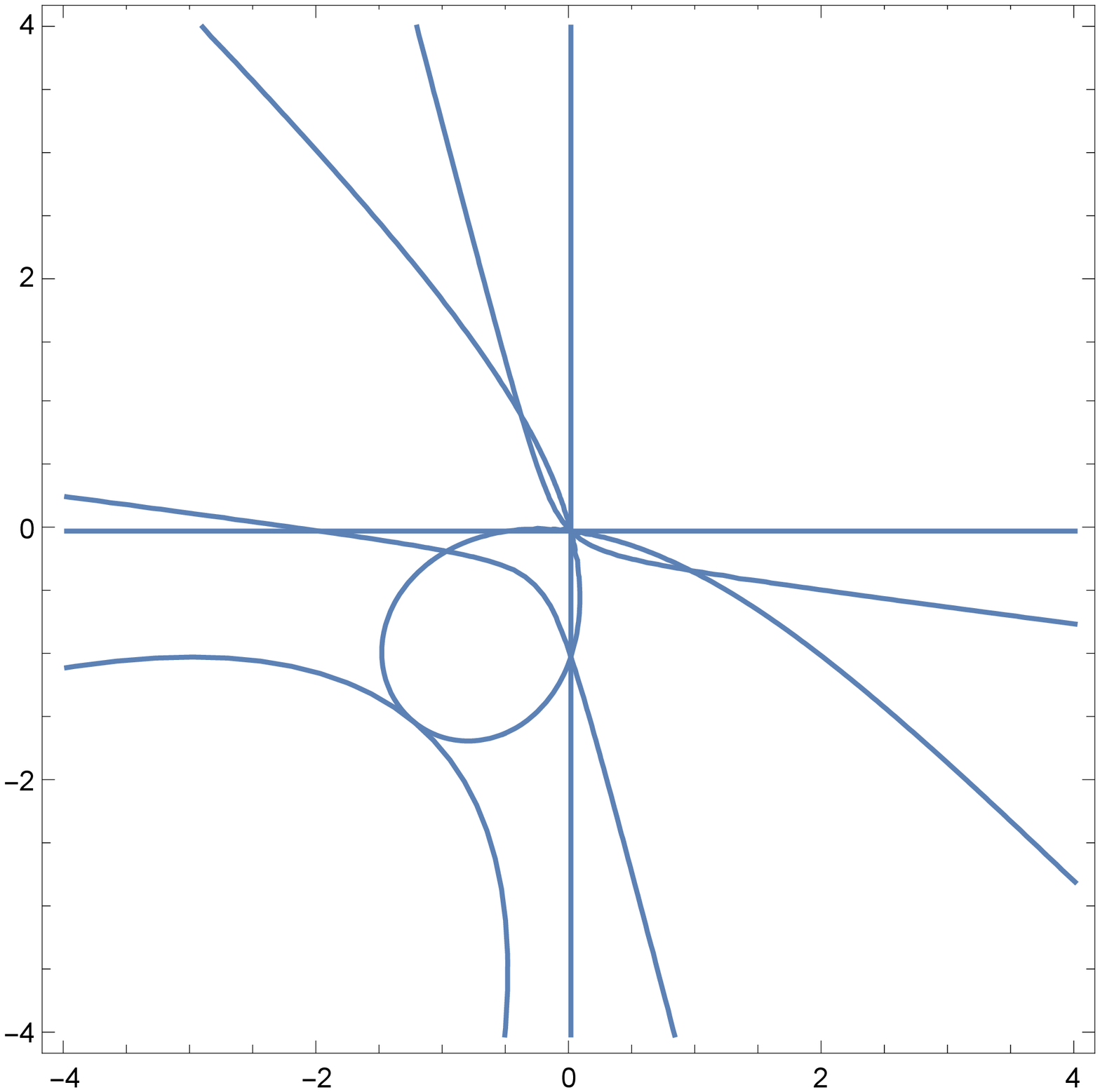}
 \\
\text{Fig. 7 : The hyperbola $\mathscr{I}(\mathscr{C}_2)$}
&
\text{Fig. 8 : Altogether}
\end{tabular}
\end{center}

\end{Exa}
\begin{Exa}
We consider a hyperbola $C_3$ (Fig. 9) defined by $3x^2-2xy+2x-2y=0.$
In this case the quartic curve $\mathscr{C}_3$ (Fig. 10) is defined by
\[
(x^2+y^2)^2+2x(x^2+y^2)+4y(x^2+y^2)+x^2+2xy+y^2=0.
\]
Then the inversion $\mathscr{I}(\mathscr{C}_3)$ (Fig. 11) of $\mathscr{C}_3$ is a parabola
defined by $
x^2+y^2+2xy+2x+4y+1=0.$
By Theorem 3.7, the quartic curve $\mathscr{C}_3$ has only one singular point (i.e. the ordinary cusp) at the origin.
The relation of these curves is depicted in Fig. 12.
\begin{center}
\begin{tabular}{cc}   
\includegraphics[width=30mm]{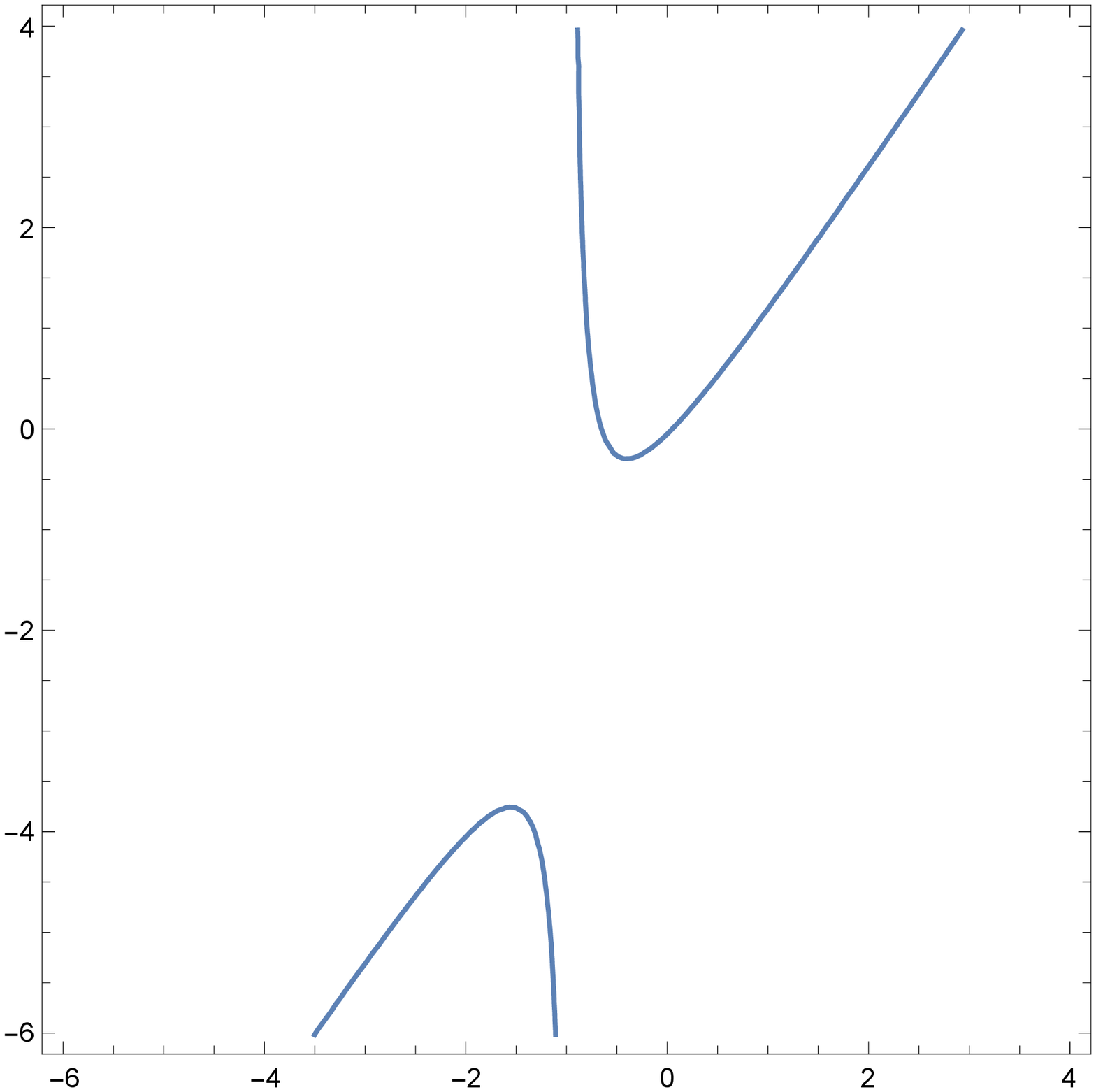}
&
\includegraphics[width=30mm]{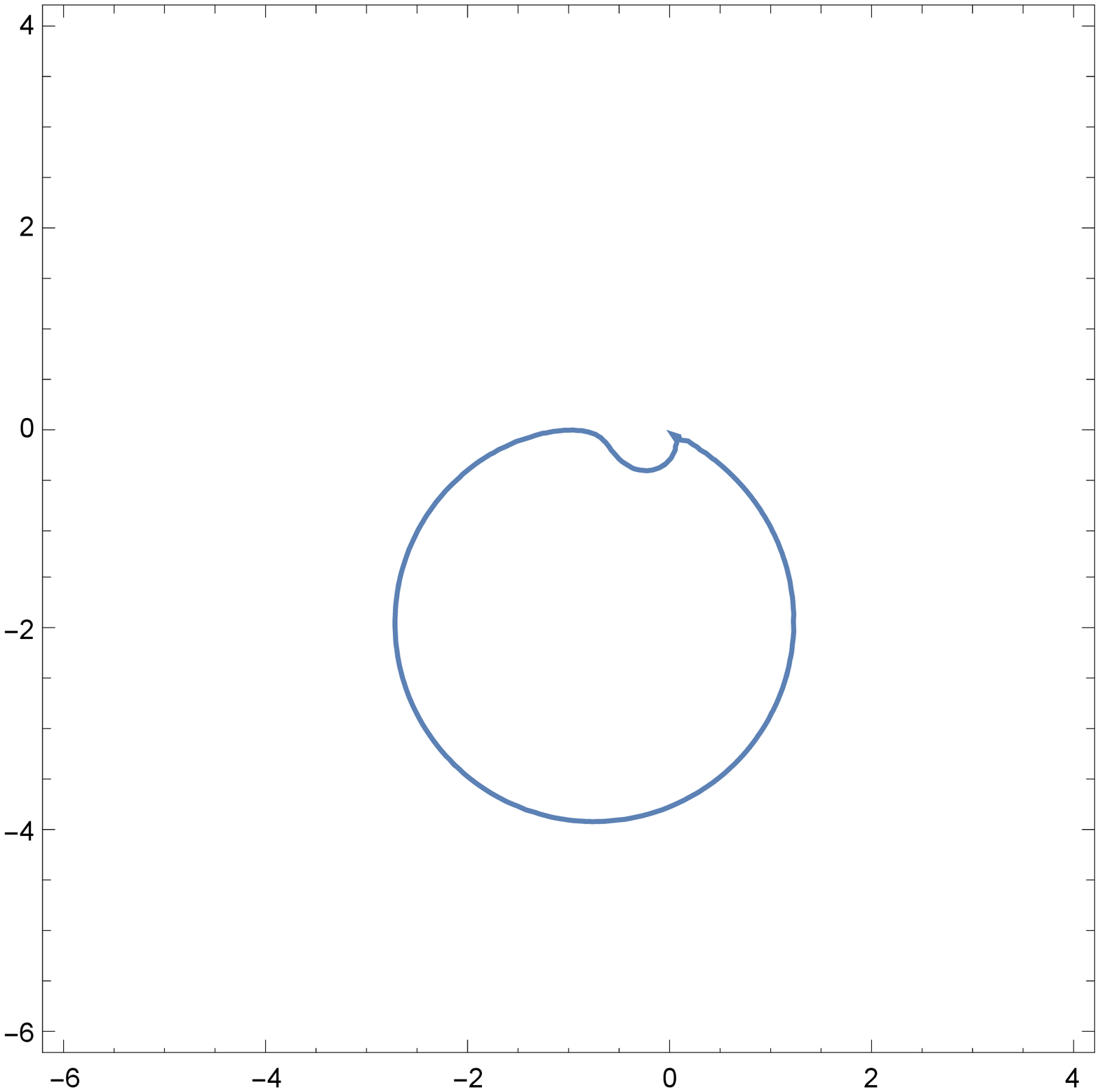}
 \\
\text{Fig. 9 : A hyperbola $C_3$}
&
\text{Fig. 10 : The quartic curve $\mathscr{C}_3$}
\end{tabular}
\end{center}

\begin{center}
\begin{tabular}{cc}   
\includegraphics[width=30mm]{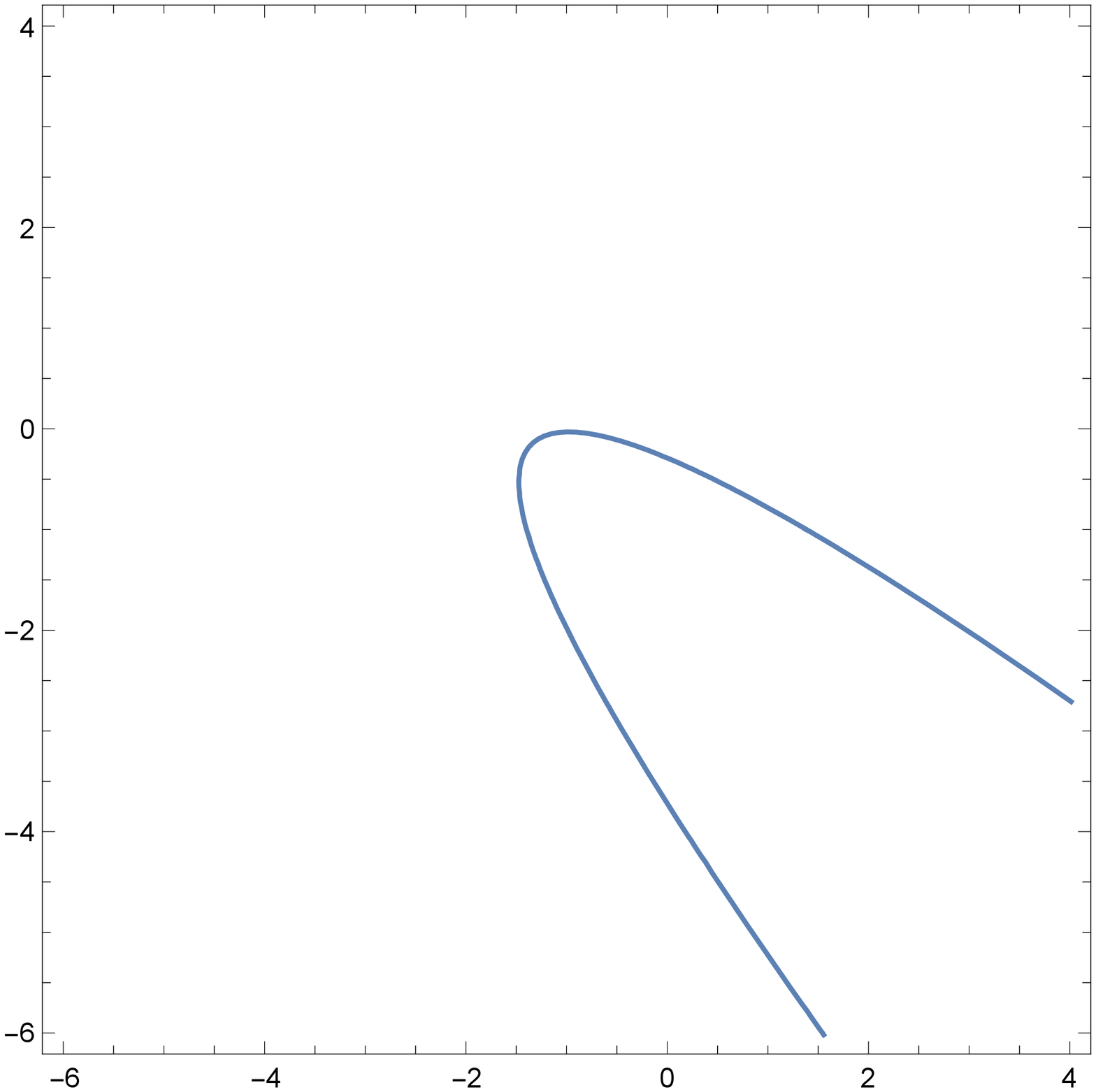}
&
\includegraphics[width=30mm]{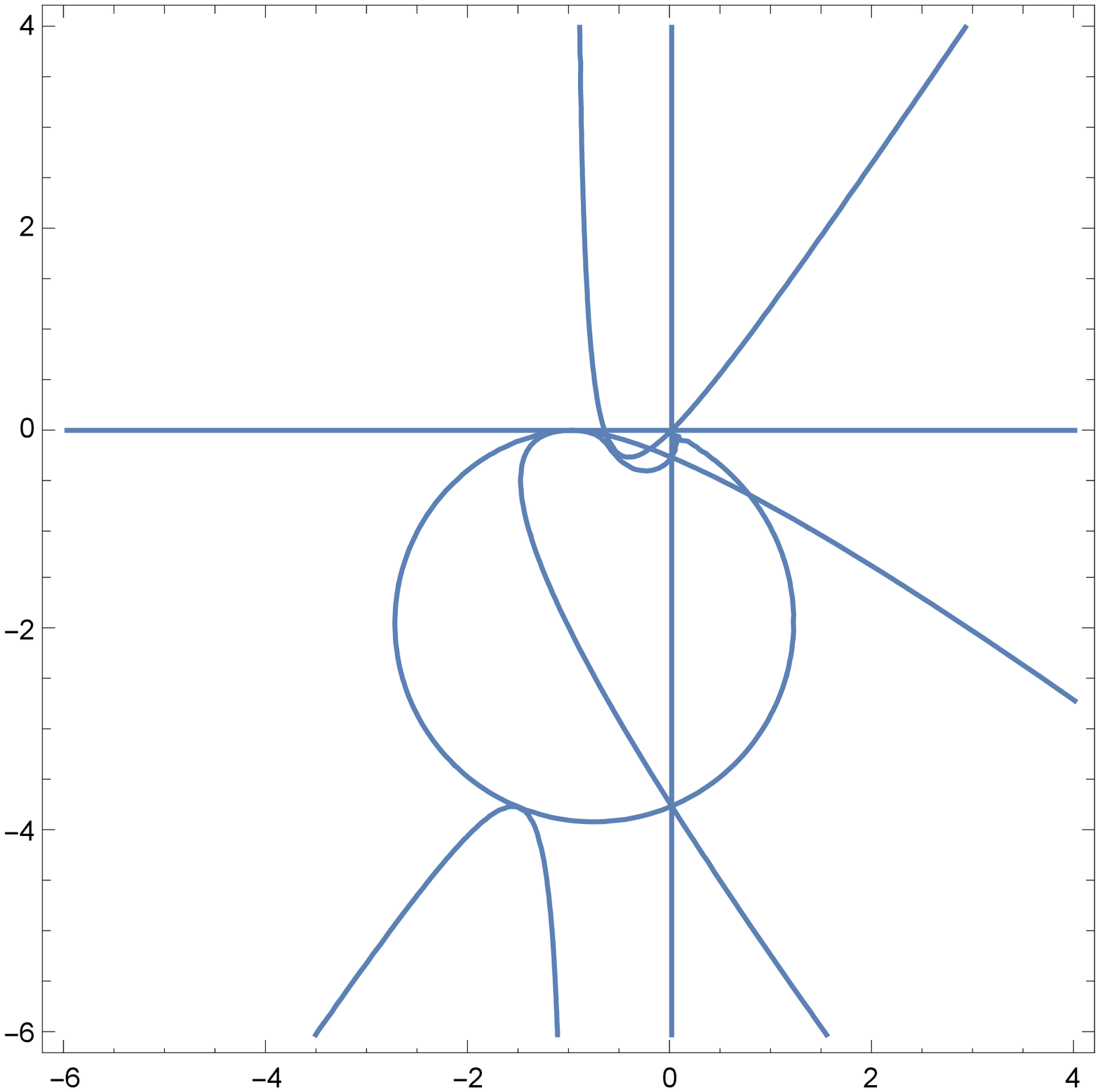}
 \\
\text{Fig. 11 : The parabola $\mathscr{I}(\mathscr{C}_3)$}
&
\text{Fig. 12 : Altogether}
\end{tabular}
\end{center}

\end{Exa}

\section{Lima\c{c}ons}
In this section we consider lima\c{c}ons. 
In the classical contexts, 
a {\it lima\c{c}on} is a plane curve parametrized by
\[
x=(r+a\cos\theta)\cos \theta,\ y=(r+a\cos\theta)\sin\theta,
\]
which is usually called a {\it lima\c{c}on of Pascal}.
We now consider a circle defined by
$(x-a)^2+(y-b)^2=r^2,$ whose parameterization is
$
\bm{\gamma}(\theta)=(r\cos\theta+a,r\sin\theta+b).
$
Then we have $\bm{t}(\theta)=(-\sin \theta, \cos\theta)$ and $\bm{n}(\theta)=(-\cos\theta,-\sin\theta).$
Therefore, $\langle \bm{\gamma}(\theta),\bm{n}(\theta)\rangle =-(r+a\cos\theta+b\sin\theta).$
Thus the pedal of the circle is
\[
{\rm Pe}_{\bm{\gamma}}(\theta)=((r+a\cos\theta+b\sin\theta)\cos\theta, (r+a\cos\theta+b\sin\theta)\sin\theta).
\]
If $(a,b)=(a,0),$ then ${\rm Pe}_{\bm{\gamma}}(\theta)$ is the lima\c{c}on of Pascal. 
On the other hand, there exists $\phi\in [0,2\pi)$ such that
$\cos\phi =a/\sqrt{a^2+b^2}, \sin\phi =b/\sqrt{a^2+b^2}.$
Then we consider the rotation
defined by 
\[
R(\phi)=\begin{pmatrix} \cos\phi &  -\sin\phi \\
\sin\phi & \cos\phi
\end{pmatrix}.
\]
It follows that
\[
R(\phi)\begin{pmatrix} r\cos\theta+a \\ r\sin\theta+b\end{pmatrix}=
\begin{pmatrix} r\cos(\theta+\phi)+\sqrt{a^2+b^2} \\ r\sin(\theta+\phi)
\end{pmatrix}. 
\]
If we set $\Theta =\theta+\phi, A=\sqrt{a^2+b^2},$
then
$
\bm{\gamma}_\phi (\Theta)=(r\cos \Theta+A, r\sin\Theta)$ is a circle with the center $(A,0).$
Therefore, the pedal of $\bm{\gamma}_\phi (\Theta)$ is the lima\c{c}on of Pascal:
\[
x=(r+A\cos\Theta)\cos \Theta,\ y=(r+A\cos\Theta)\sin\Theta.
\]
By definition, we have $R(\phi){\rm Pe}_{\bm{\gamma}}={\rm Pe}_{\bm{\gamma}_\phi}.$
This means that the pedal of $\bm{\gamma}$ is given by the rotation of the lima\c{c}on of Pascal.
Therefore, a (generalized) {\it lima\c{c}on} is defined to be
\[
((r+a\cos\theta+b\sin\theta)\cos\theta, (r+a\cos\theta+b\sin\theta)\sin\theta),
\]
which is the pedal of a circle.
By the arguments in \S 3, the implicit equation of the lima\c{c}on has the following form:
\[
(x^2+y^2)^2-2(ax+by)(x^2+y^2)+(a^2-r^2)x^2+(r^2-b^2)y^2-2abxy+2(ax+by)-1=0
\]
as the pedal of the circle $(x-a)^2+(y-b)^2=r^2.$
We remark that lima\c{c}ons are always quartic curves and fronts. 
 We have the following theorem as a corollary of Theorems 3.1 and 3.7.
\begin{Th} Let $\mathscr{C}$ be an algebraic curve in $\R^2$ such that $C=\mathcal{P}r_{\mathscr{C}^*}$ is a circle.
Then we have the following assertions\/{\rm :}
\par\noindent
{\rm (1)} The algebraic curve $\mathscr{C}$ is a lima\c{c}on such that 
$\mathscr{I}(\mathscr{C})$ is an irreducible quadratic curve.
Moreover, if the equation of $C$ is 
$
(x-a)^2+(y-b)^2=r^2,
$
then the equation of $\mathscr{C}$ is 
\[
(x^2+y^2)^2-2(ax+by)(x^2+y^2)+(a^2-r^2)x^2+(b^2-r^2)y^2+2abxy=0,
\]
and the equation of the inversion $\mathscr{I}(\mathscr{C})$ of $\mathscr{C}$ is
\[
(r^2-a^2)x^2+(r^2-b^2)y^2-2abxy+2(ax+by)-1=0.
\]
\par\noindent
{\rm (2)} If $a^2+b^2-r^2=0,$ then $\mathscr{I}(\mathscr{C})$ is a parabola and
the lima\c{c}on $\mathscr{C}$ has  an ordinary cusp at the origin.
\par\noindent
{\rm (3)} If $a^2+b^2-r^2<0,$ then $\mathscr{I}(\mathscr{C})$ is an ellipse and
the lima\c{c}on $\mathscr{C}$ has an isolated point at the origin and other points are non-singular.
\par\noindent
{\rm (4)} If $a^2+b^2-r^2>0,$ then $\mathscr{I}(\mathscr{C})$ is a hyperbola and
the lima\c{c}on $\mathscr{C}$ has a node at the origin and any point is non-singular as a front.
\end{Th}
\demo
By assertion (1) of Theorem 3.7, assertion (1) holds.
By the equation of $\mathscr{I}(\mathscr{C}),$ we have
$\Delta _0=r^2(r^2-a^2-b^2)$ and $\Delta =-r^4.$
By Theorems 3.1 and 3.7, the other assertions hold.
This completes the proof.
\enD
\begin{Exa} We give some examples of lima\c{c}ons.
\par\noindent
(1) We consider a circle $C_4$ (Fig. 13) defined by $(x-3)^2+y^2=9.$
In this case the quartic curve $\mathscr{C}_4$ (i.e. the lima\c{c}on; Fig. 14) is defined by
\[
(x^2+y^2)^2-6x(x^2+y^2)-9y^2=0.
\]
Then the inversion $\mathscr{I}(\mathscr{C}_4)$ (Fig. 15) of $\mathscr{C}_4$ is a a parabola
defined by $
9y^2+6x-1=0.$
The relation of these curves is depicted in Fig. 16.
\begin{center}
\begin{tabular}{cc}   
\includegraphics[width=35mm]{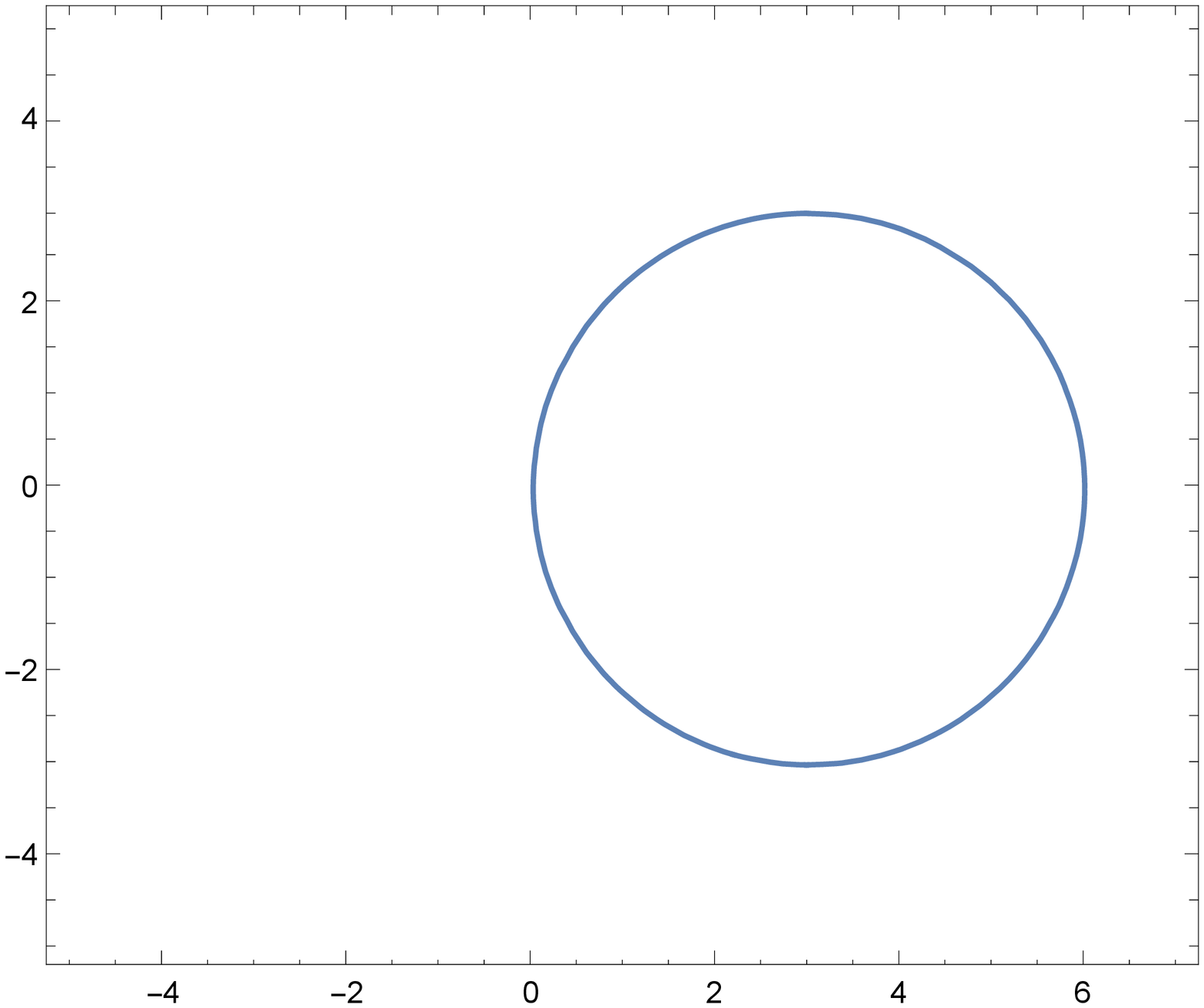}
&
\includegraphics[width=35mm]{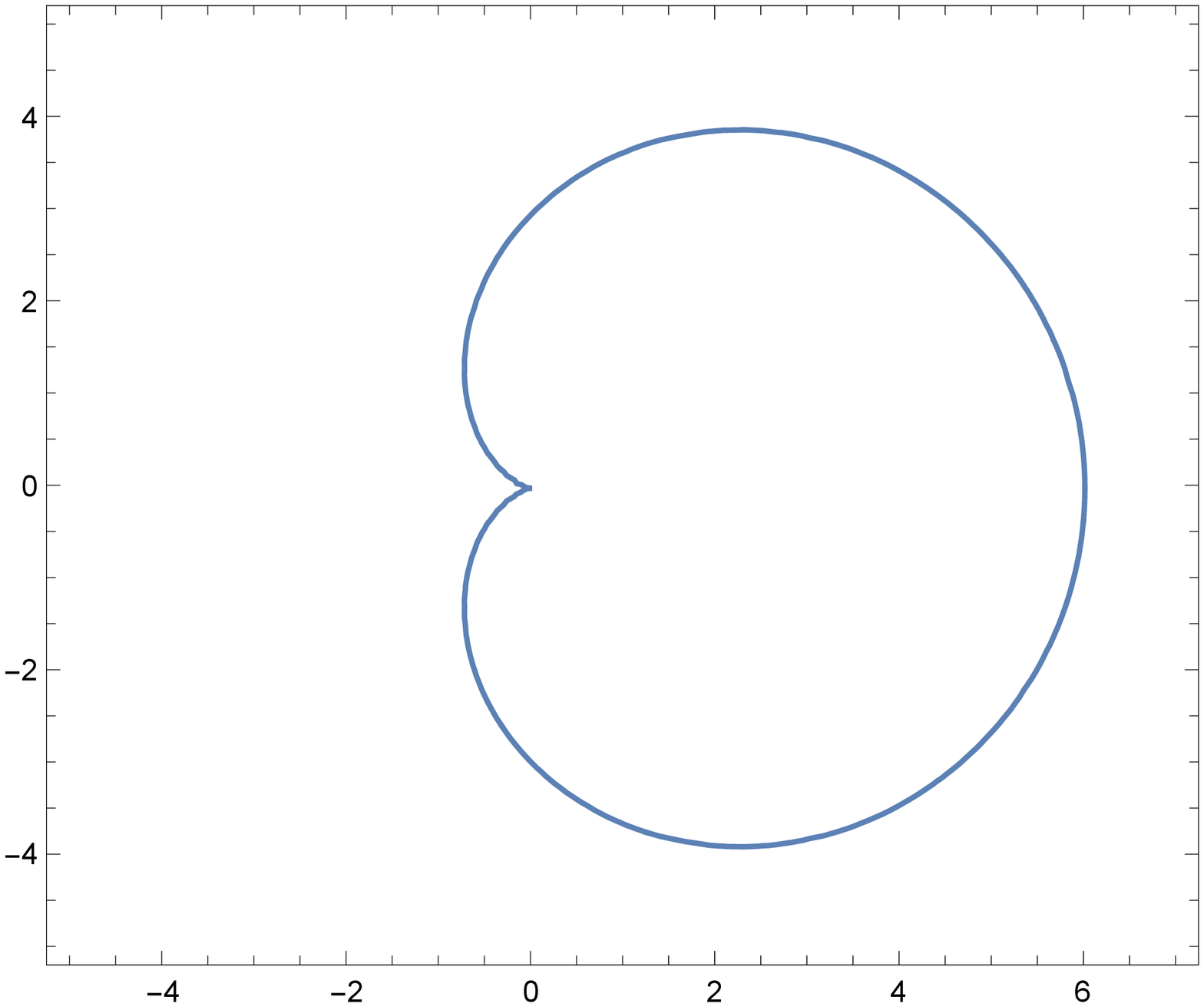}
 \\
\text{Fig. 13 : The circle $C_4$}
&
\text{Fig. 14 : The lima\c{c}on $\mathscr{C}_4$}
\end{tabular}
\end{center}

\begin{center}
\begin{tabular}{cc}   
\includegraphics[width=35mm]{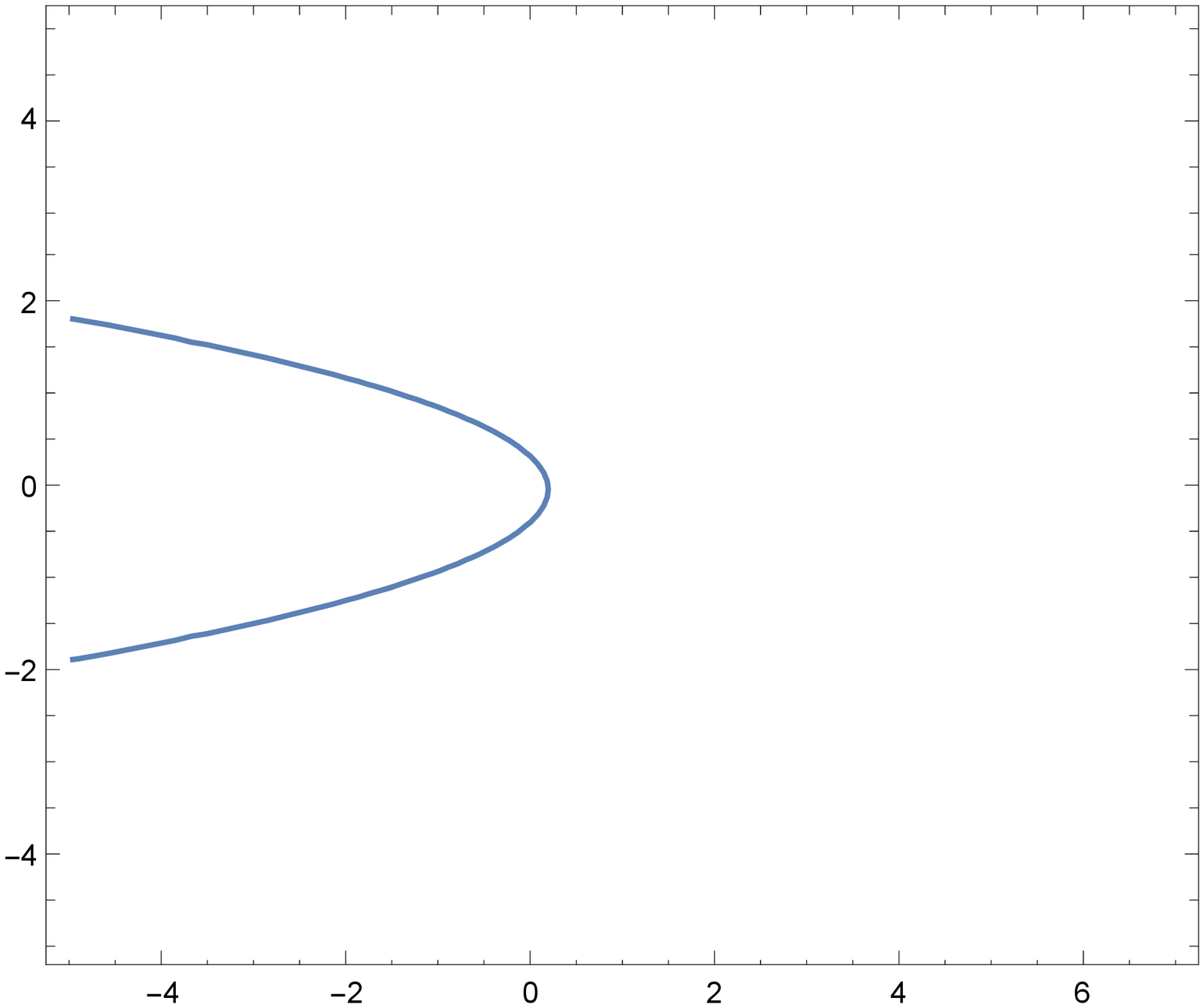}
&
\includegraphics[width=35mm]{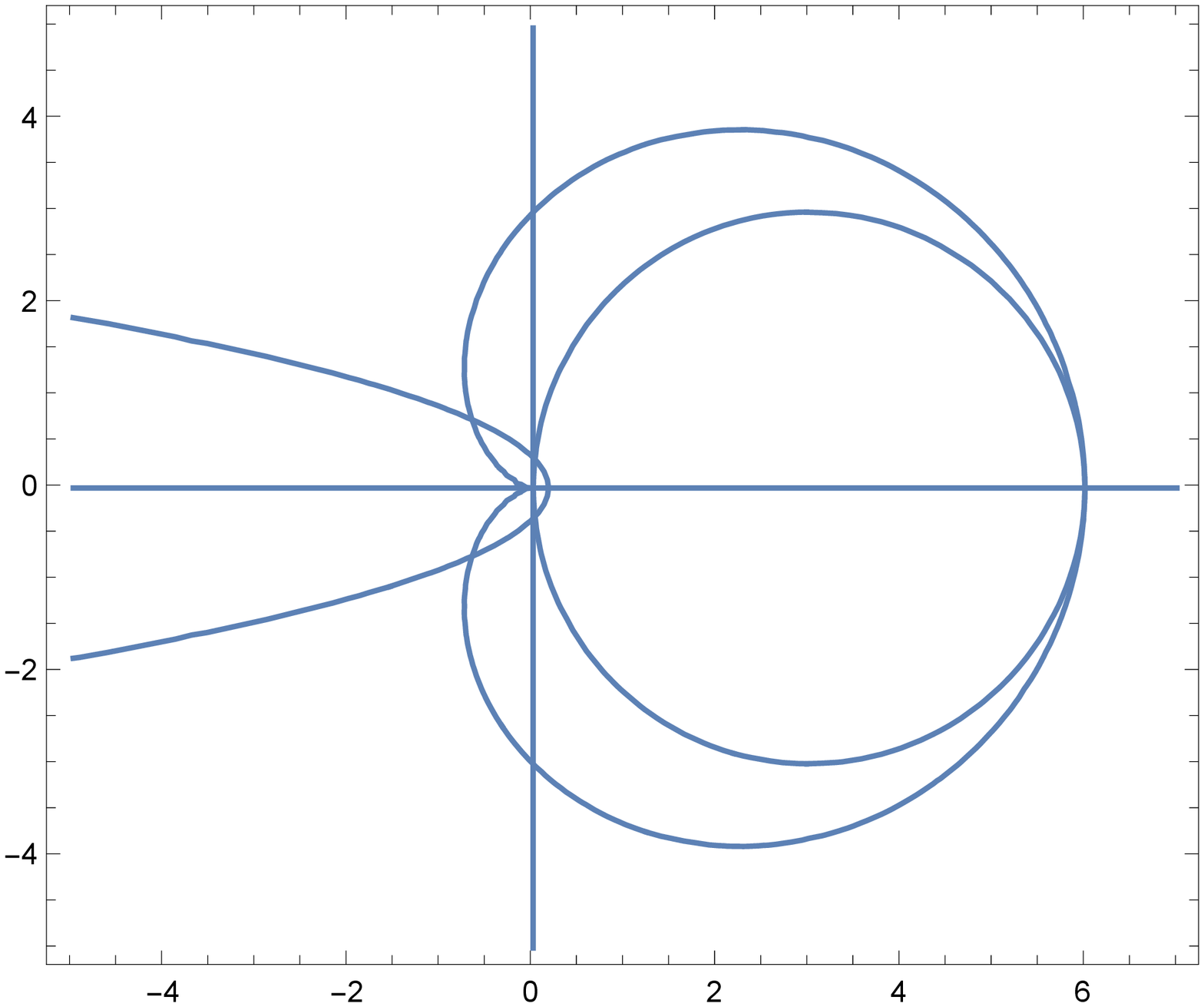}
 \\
\text{Fig. 15 : The parabola $\mathscr{I}(\mathscr{C}_4)$}
&
\text{Fig. 16 : Altogether}
\end{tabular}
\end{center}
\par\noindent
(2) We consider a circle $C_5$ (Fig. 17) defined by $(x-1)^2+(y-2)^2=7.$
In this case the quartic curve $\mathscr{C}_5$ (i.e. the lima\c{c}on; Fig. 18)) is defined by
\[
(x^2+y^2)^2-(2x+4y)(x^2+y^2)-6x^2-3y^2+4xy=0.
\]
Since $a^2+b^2-r^2=-2<0,$ the inversion $\mathscr{I}(\mathscr{C}_5)$ (Fig. 19) of $\mathscr{C}_5$ is a an ellipse
defined by 
\[
6x^2+3y^2-4xy+2x+4y-1=0.
\]
The relation of these curves is depicted in Fig. 20.
\begin{center}
\begin{tabular}{cc}   
\includegraphics[width=30mm]{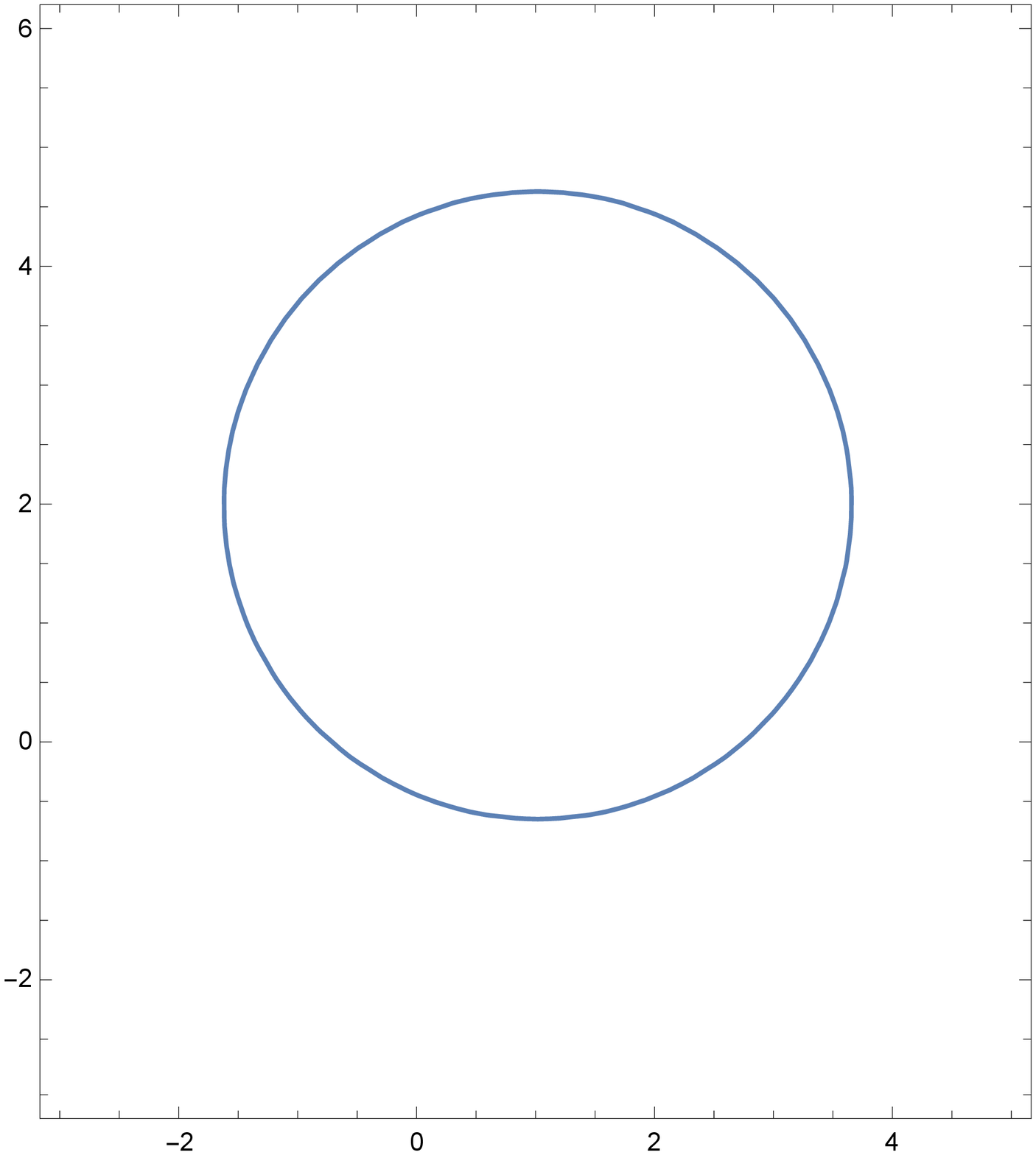}
&
\includegraphics[width=30mm]{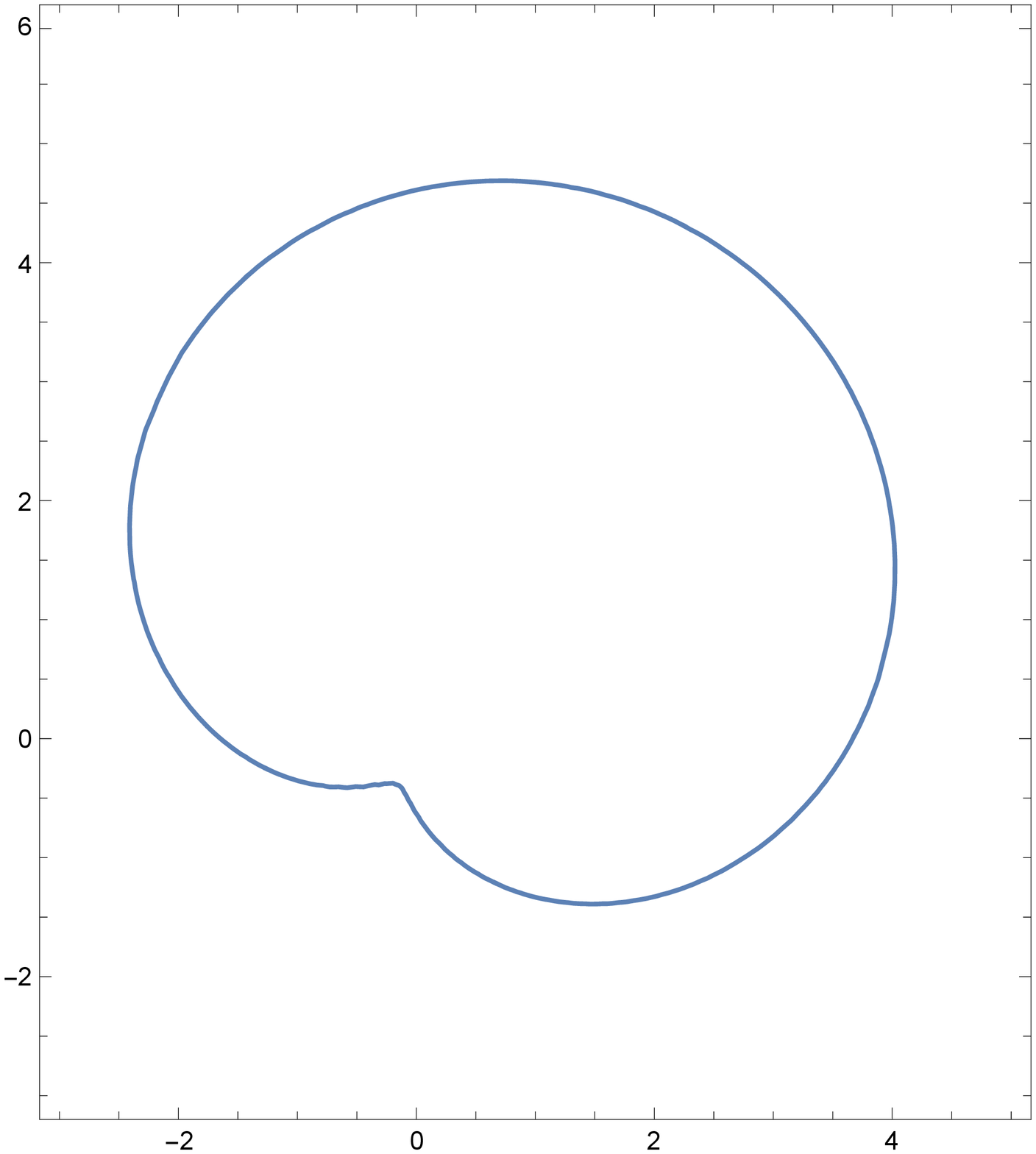}
 \\
\text{Fig. 17 : The circle $C_5$}
&
\text{Fig. 18 : The lima\c{c}on $\mathscr{C}_5$}
\end{tabular}
\end{center}

\begin{center}
\begin{tabular}{cc}   
\includegraphics[width=30mm]{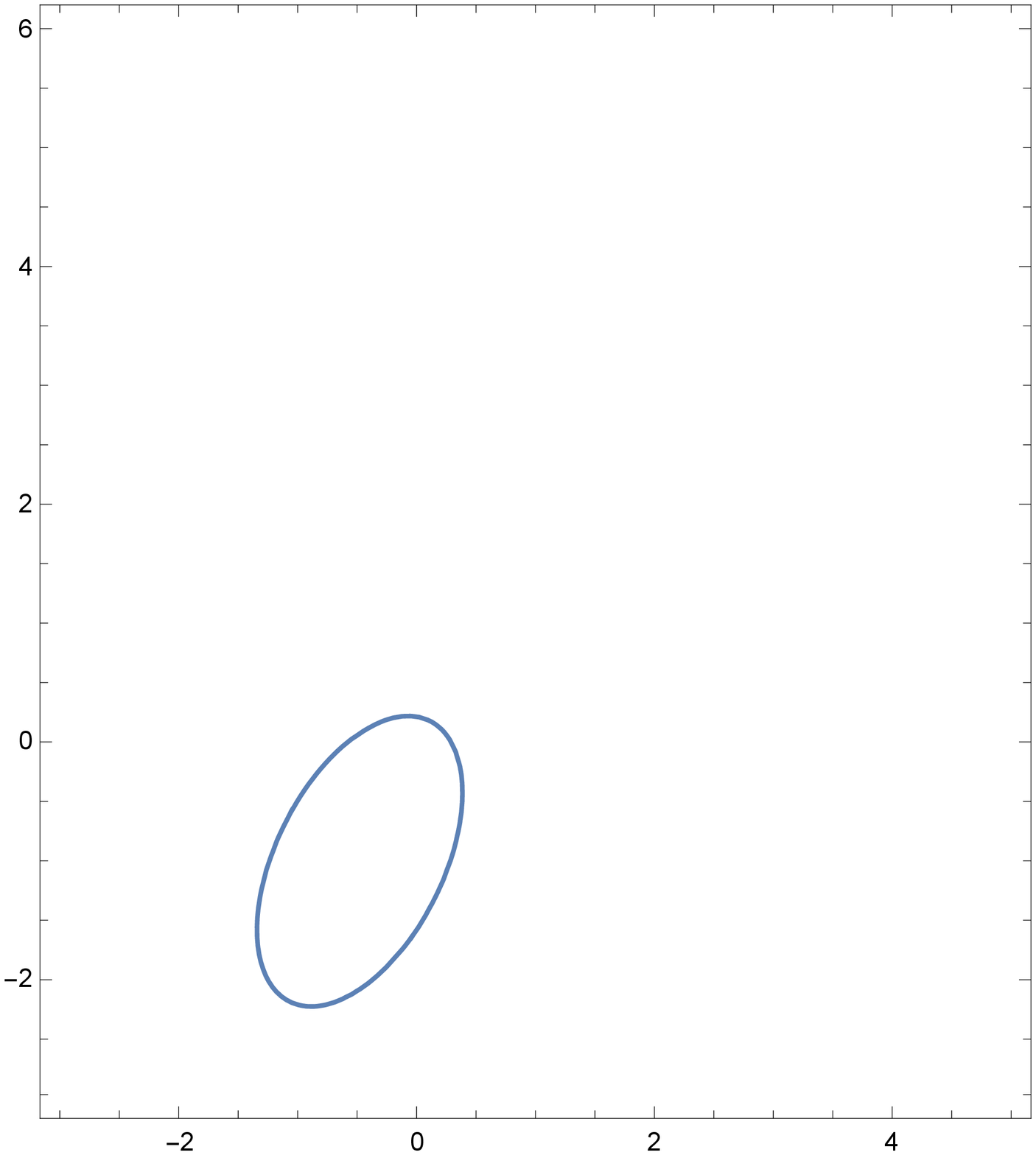}
&
\includegraphics[width=30mm]{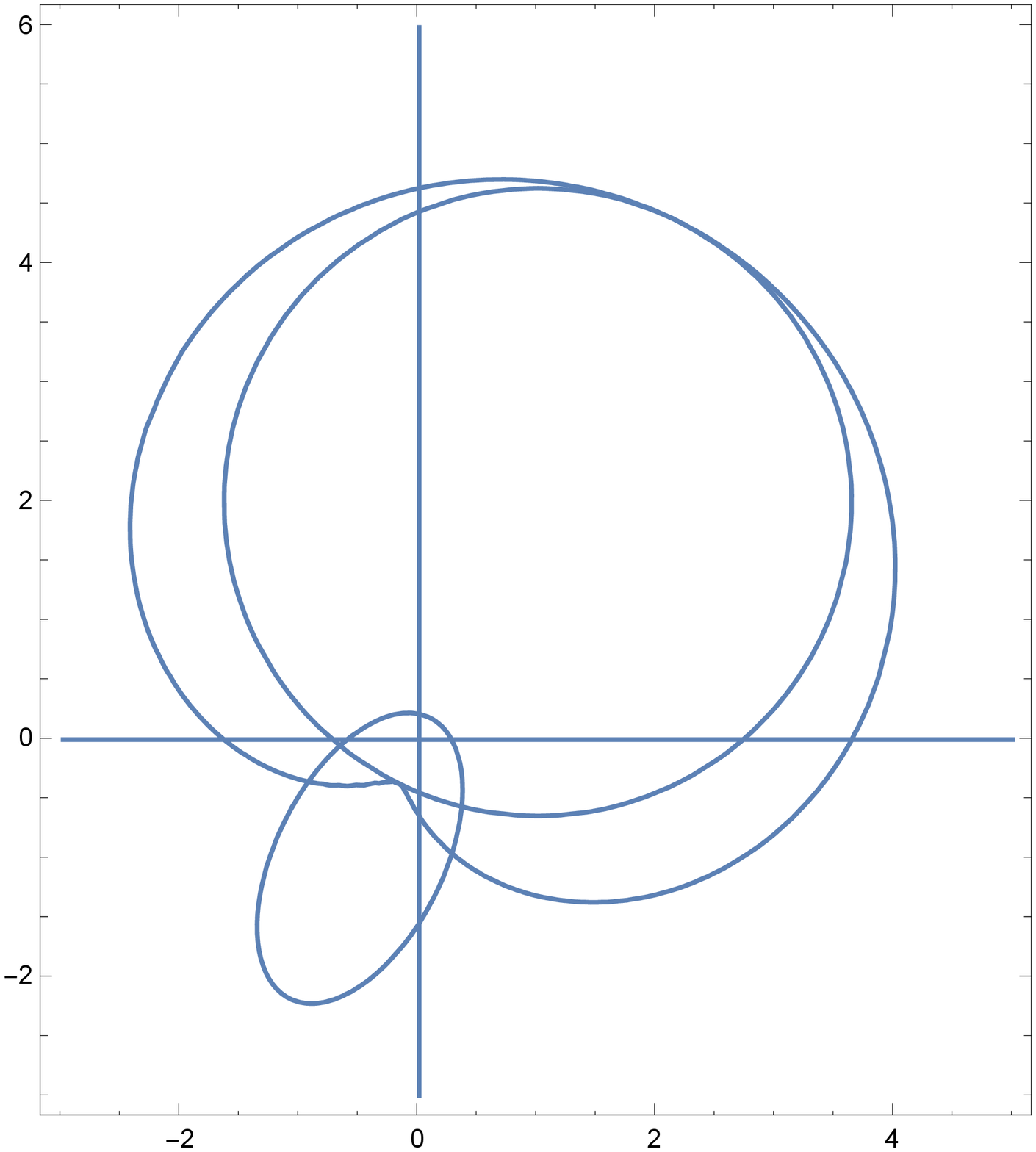}
 \\
\text{Fig. 19 : The parabola $\mathscr{I}(\mathscr{C}_5)$}
&
\text{Fig. 20 : Altogether}
\end{tabular}
\end{center}
\par\noindent
 (3) We consider a circle $C_6$ (Fig. 21) defined by $(x-1)^2+(y-2)^2=1.$
In this case the quartic curve $\mathscr{C}_6$ (i.e. the lima\c{c}on; Fig. 22) is defined by
\[
(x^2+y^2)^2-(2x+4y)(x^2+y^2)+3y^2+4xy=0.
\]
Since $a^2+b^2-r^2=4>0,$ the inversion $\mathscr{I}(\mathscr{C}_6)$ (Fig. 23) of $\mathscr{C}_6$ is a a hyperbola
defined by 
\[
3y^2+4xy-2x-4y+1=0.
\]
The relation of these curves is depicted in Fig. 24.
\begin{center}
\begin{tabular}{cc}   
\includegraphics[width=35mm]{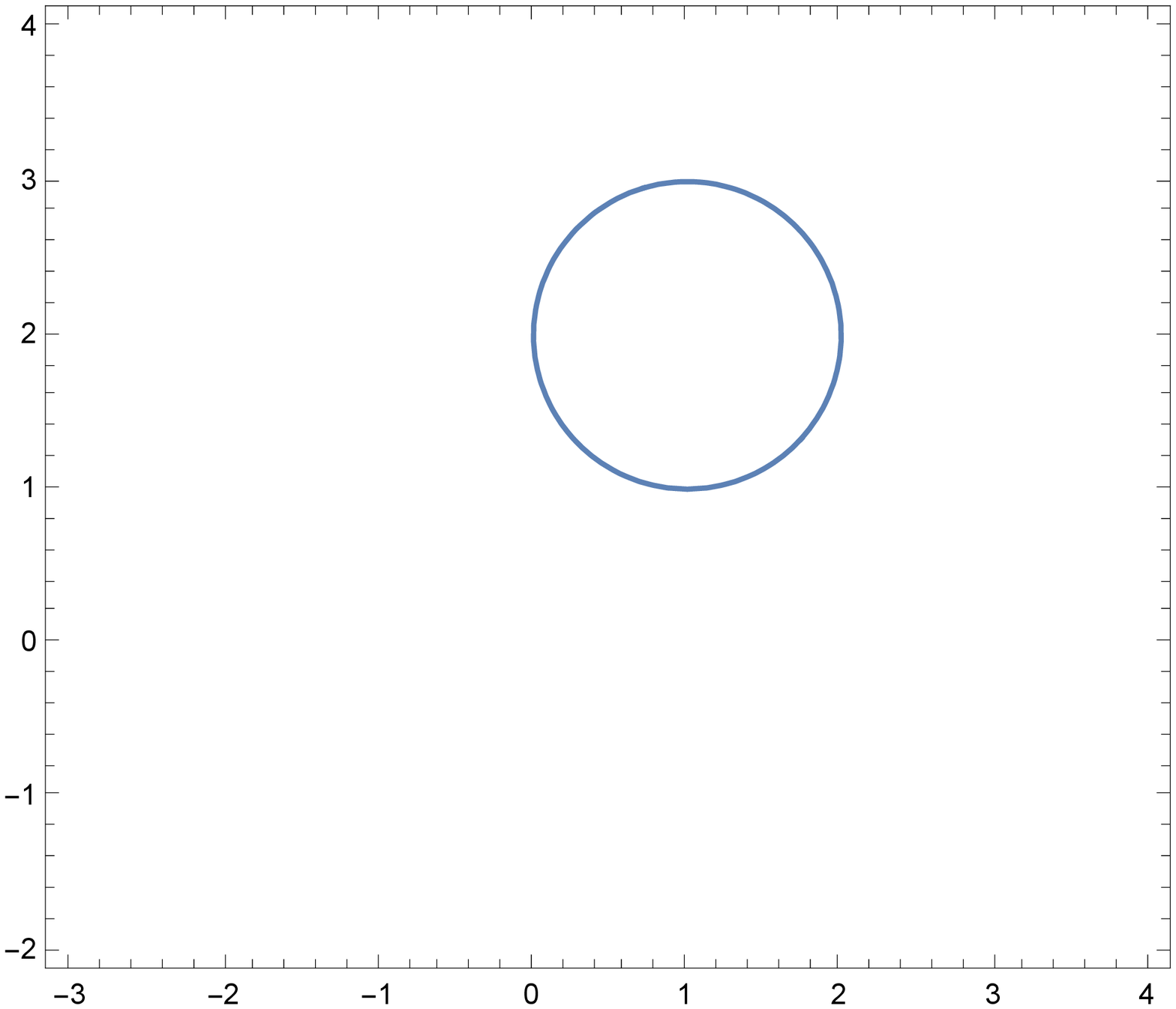}
&
\includegraphics[width=35mm]{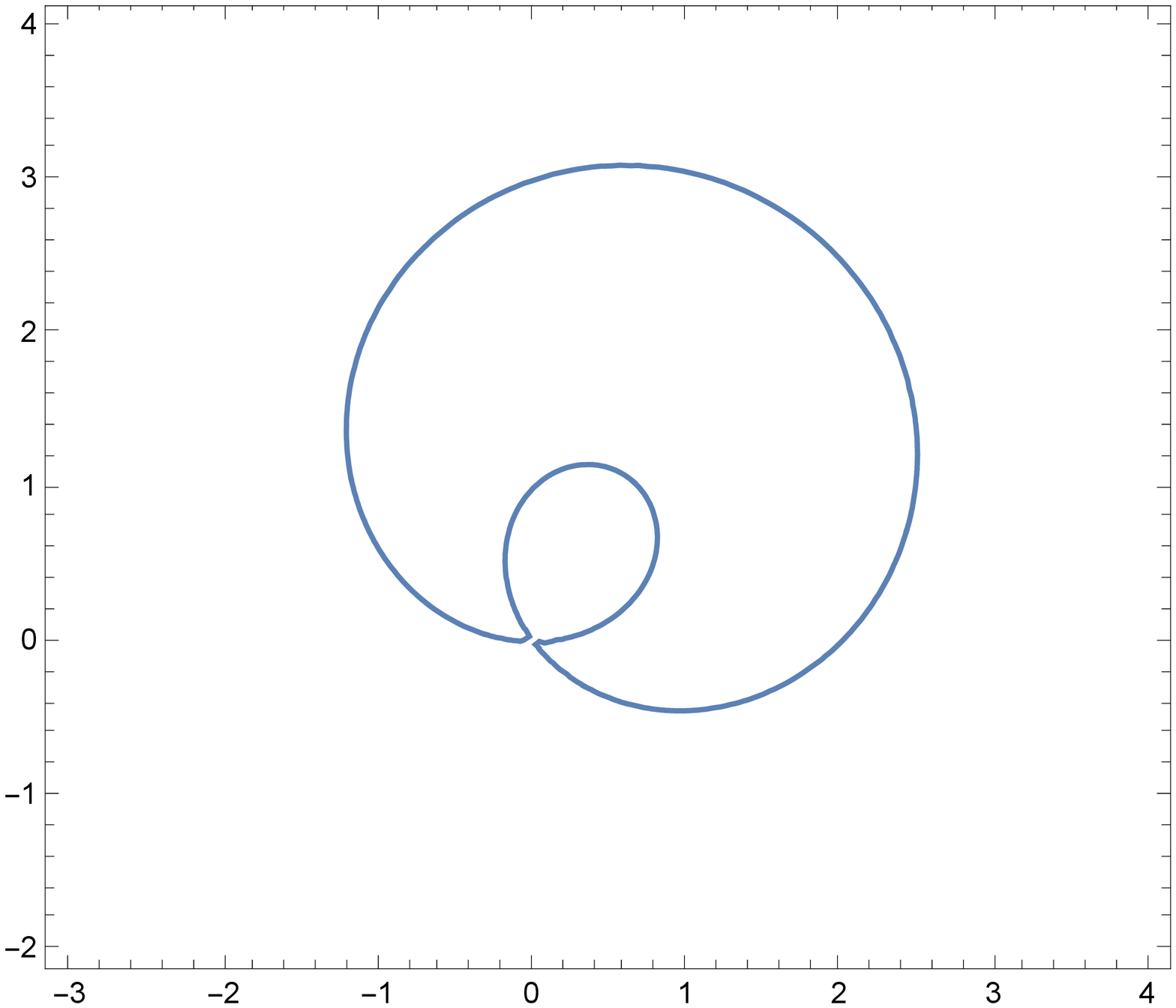}
 \\
\text{Fig. 21 : The circle $C_6$}
&
\text{Fig. 22 : The lima\c{c}on $\mathscr{C}_6$}
\end{tabular}
\end{center}

\begin{center}
\begin{tabular}{cc}   
\includegraphics[width=35mm]{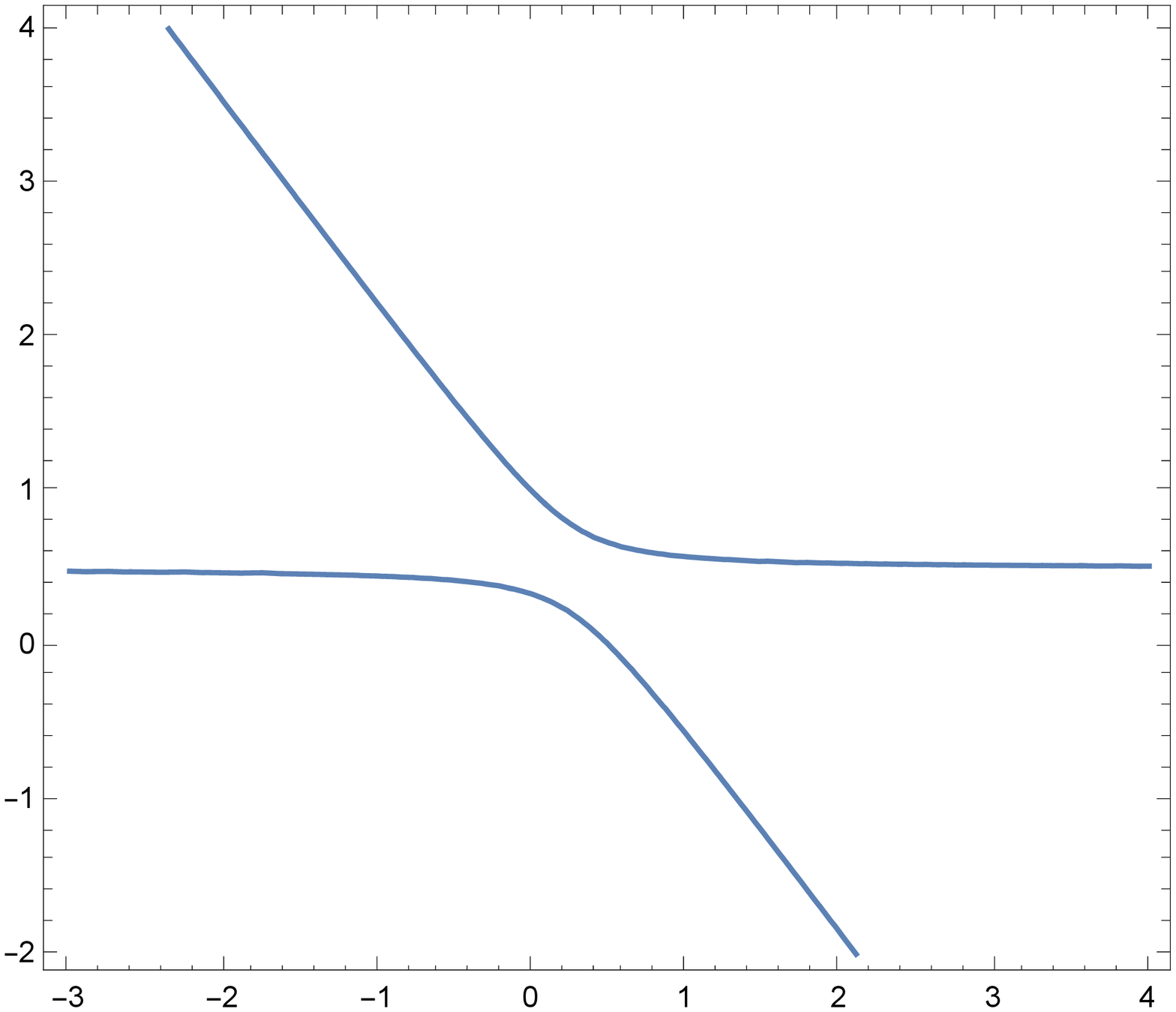}
&
\includegraphics[width=35mm]{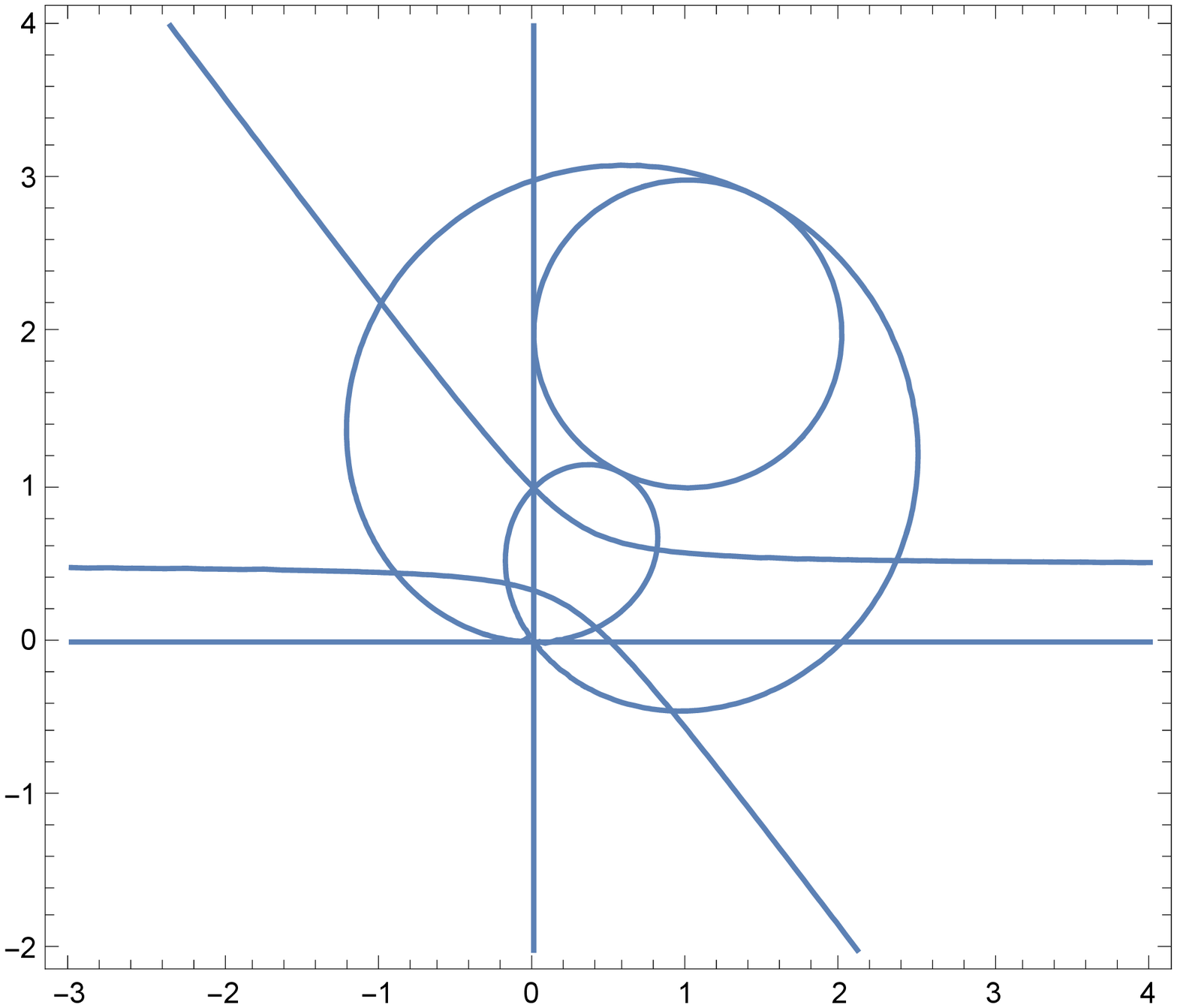}
 \\
\text{Fig. 23 : The parabola $\mathscr{I}(\mathscr{C}_6)$}
&
\text{Fig. 24 : Altogether}
\end{tabular}
\end{center}

\end{Exa}

\begin{flushright}
\begin{tabular}{l}
Shyuichi Izumiya\\
Department of Mathematics, \\
Hokkaido University,\\
Sapporo 060-0810, Japan \hspace*{26mm}\\
{\tt izumiya@math.sci.hokudai.ac.jp}
\end{tabular}
\end{flushright}
\begin{flushright}
\begin{tabular}{l}
Nobuko Takeuchi\\
Department of Mathematics,\\
Tokyo Gakugei University,\\
Koganei, Tokyo, 184-8501, Japan \hspace*{12mm}\\
{\tt nobuko@u-gakugei.ac.jp}\\
\end{tabular}
\end{flushright}

\end{document}